\documentclass[final,3p,times]{elsarticle}

\usepackage{amssymb}
\usepackage{amsmath}
\usepackage{amsthm}
\usepackage{here}
\usepackage{booktabs} 
\usepackage{tabularx} 
\newtheorem{prop}{Proposition}[section]
\newtheorem{thm}{Theorem}[section]

\newtheorem{rem}{Remark}[section]
\newtheorem{cor}{Corollary}[section]

\newtheorem{dfn}{Definition}[section]

\newtheorem{exm}{Example}[section]
\newcommand{\Eqref}[1]{Eq.~(\ref{#1})}

\newcommand{\Section}[1]{Sec.~\ref{#1}}
\newcommand{\Figure}[1]{Fig.~\ref{#1}}
\newcommand{\Theorem}[1]{Theorem~\ref{#1}}
\newcommand{\Corollary}[1]{Corollary~\ref{#1}}
\newcommand{\Proposition}[1]{Proposition~\ref{#1}}

\newcommand{\Table}[1]{Table~\ref{#1}}
\newcommand{\discFunc}[2]{{#1}^{(#2)}}
\newcommand{\discreteFunc}[4]{{#1}_{#2,#3}^{(#4)}}

\newcommand{\discreteScheme}[2]{S_{#1}^{#2}}

\newcommand{\numBulk}[0]{N}
\newcommand{\numBdd}[0]{N^\partial}
\newcommand{\thresTildePM}[0]{\widetilde{\alpha}_s^\pm}
\newcommand{\thresTildeP}[0]{\widetilde{\alpha}_s^+}
\newcommand{\thresTildeM}[0]{\widetilde{\alpha}_s^-}
\newcommand{\thresP}[0]{\alpha_s^+}
\newcommand{\thresM}[0]{\alpha_s^-}
\newcommand{\thresPM}[0]{\alpha_s^\pm}
\newcommand{\intval}[2]{\left(#1,#2\right)}

\newcommand{\intvalClosed}[2]{[#1,#2]}

\newcommand{\trace}[1]{\operatorname{Tr}{\left(#1\right)}}

\newcommand{\closure}[1]{\overline{#1}}
\newcommand{\detK}[0]{\Delta}
\newcommand{\matK}[0]{\left[k\right]}
\newcommand{\mat}[1]{#1}
\newcommand{\pOmega}[0]{\partial\Omega}

\makeatletter
\def\Hline{%
\noalign{\ifnum0=`}\fi\hrule \@height 1.5pt \futurelet
\reserved@a\@xhline}
\makeatother

\journal{Computers \& Mathematics with Applications}

\begin{document}

\begin{frontmatter}



\title{A hyperbolic finite difference scheme for anisotropic diffusion equations: preserving the discrete maximum principle}


\author[label1]{Tokuhiro Eto\corref{cor1}}
\ead{eto@math.univ-lyon1.fr}
\author[label2]{Rei Kawashima}

\cortext[cor1]{Corresponding author.}
\ead{reik@shibaura-it.ac.jp}
\address[label1]{Universit\'{e} Claude Bernard Lyon 1, CNRS, Centrale Lyon, INSA Lyon, Universit\'{e} Jean Monnet, ICJ UMR5208, 69622 Villeurbanne, France}
\address[label2]{Department of Electrical Engineering, Shibaura Institute of Technology, 3-8-5 Toyosu, Koto-ku, Tokyo 135-8548, Japan}

\begin{abstract}
	A hyperbolic system approach is proposed for robust computation of anisotropic diffusion equations that appear in quasineutral plasmas.
	Though the approach exhibits merits of high extensibility and accurate flux computation, the monotonicity of the scheme for anisotropic diffusion cases has not been understood.
	In this study, the discrete maximum principle (DMP) of the hyperbolic system approach is analyzed and tested in various anisotropic diffusion cases.
	A mathematical analysis is conducted to obtain an optimal condition of an arbitrary parameter to guarantee the DMP, and numerical experiments reveal an adoptive selection of the parameter for DMP-preserving results.
	It is confirmed that, with an appropriate preconditioning matrix and parameter choice, the hyperbolic system approach preserves the DMP even with a linear discretization.
\end{abstract}

\end{frontmatter}



\section{Introduction}\label{sec:intro}
	The plasma sources that employ applied electric and magnetic fields are used in a variety of applications, inclusding electric propulsion, semiconductor manufacturing, and nuclear fusion.
	In these devices,  the magnetic field plays an essential role in magnetic confinement and in controlling the motion of charged particles.
	Both the physical model and numerical method must accurately capture the effects of the applied magnetic field in the computer-aided-engineering (CAE) processes to develop such plasma sources.
	In many plasma sources, a partially magnetized plasma is generated: electrons are strongly magnetized, whereas ions remain essentially unmagnetized to enable smooth extraction.
	A low-temperature ($<100$ eV) plasma is generated by introducing a operating gas into the discharge chamber.
	The resulting plasma flows in industrial plasma sources are typically partially ionized, partially magnetized, low-temperature, and collisional.
	In this regime, the electron flow becomes a diffusion-dominated fluid with strong anisotropy owing to magnetic confinement.
	The governing equation for magnetized electron flow therefore takes the form of an  anisotropic diffusion equation, which is the focus of this study.

	An anisotropic diffusion equation is commonly defined as follows:
	\begin{align}\label{eq:ade}
		\begin{split}
			&\nabla \cdot \vec{\Gamma} =0, \\
			&\vec{\Gamma}= -\left[k\right]\nabla\phi,
		\end{split}
	\end{align}
	where $\left[k\right]$ is the diffusion-coefficient matrix which is $2\times2$ for two-dimensional problems.
	The quantity $\phi$ is the main variable of the equation and corresponds to the space potential in the case of magnetized electron flow.
	The vector $\vec{\Gamma}$ denotes the flux arising from the gradients of the main variable.
	The coefficient matrix in \Eqref{eq:ade} is written as
	\begin{equation}
		\left[k\right]=\left[\begin{array}{cc}
				k_x & k_{\rm c} \\
				k_{\rm c} & k_y
			\end{array}\right]
		=\Theta^{-1}
			\left[\begin{array}{cc}
				k_{||} & \ \\
				\ & k_{\perp}
			\end{array}\right]\Theta,
		\label{eq:kmat}
	\end{equation}
	where $\Theta$ is the rotation matrix with respect to the angle between the magnetic-field line and the computational grid.
	The strength of the anisotropy is the ratio of the diffusion coefficients parallel and perpendicular to the field, $k_{||}/k_{\perp}$.
	The main difference of this problem from standard isotopic diffusion is that the diffusion-coefficient matrix contains off-diagonal elements, $k_{\rm c}$, which may cause numerical instabilities during iterative calculation or spurious oscillations in the results.
	It is known that the diagonal dominance (or positivity) is lost in a numerical scheme based on a standard central differencing for the cross-diffusion terms \cite{kawashima2015}.
	Consequently, many practical simulations of magnetized plasmas---such as tokamak plasmas and Hall thrusters---avoid the treatment of off-diagonal terms by employing a magnetic-field-aligned mesh (MFAM) \cite{leddy2017,mikellides2011}.
	However, constructing an MFAM is difficult for complex magnetic-field geometries, such as cusped magnetic fields or fields distorted by a plasma-induced currents.
	To retain flexibility with respect to magnetic-field geometry, a numerical method that handles the off-diagonal terms in the anisotropic diffusion equation is required.

	The hyperbolic-equation-system approach has been developed for robust and accurate computations of anisotropic-diffusion problems of magnetized electrons.
	In this approach, the second-order diffusion equation is converted into an equivalent system of first-order partial-differential equations.
	The resulting hyperbolic system can be computed robustly using numerical schemes developed in computational fluid dynamics (CFD) for compressible fluids.
	A similar concept was first proposed for isotropic diffusion and advection-diffusion equations \cite{nishikawa2007,nishikawa2010}.
	Kawashima et al. developed a hyperbolic system approach suitable for anisotropic diffusion of the magnetized, and demonstrated that it stably computes the cross-diffusion terms \cite{kawashima2015}.
	Furthermore, a series of studies has proven the following characteristics of the hyperbolic system approach.
	1) High extensibility: the approach can be extended to equation sets that include several advection and diffusion terms.
	In fact, this approach has been extended to the governing equations of magnetized electrons (mass, momentum, and internal energy) in a quasineutral plasma \cite{kawashima2016}, yielding six first-order equations in two dimensions.
	Its applicability to the three-dimensional Navier--Stokes equation has also been confirmed \cite{nakashima2016}, involving 20 first-order equations.
	2) High accuracy in flux (gradient) calculation: numerical fluxes are strictly conserved in the finite volume framework.
	Furthermore, the gradients of the variables are computed with the same order of accuracy as the primary variable \cite{nishikawa2010}.
	This property is particularly important in magnetized plasmas, where cross-field fluxes directly relate to magnetic-confinement performance.
	Owing to these advantages the hyperbolic system approach is expected to be beneficial for a practical simulation of electric propulsion systems \cite{kawashima2019}.

\begin{figure}[t]
	   	\begin{center}
   			\includegraphics[width=60mm]{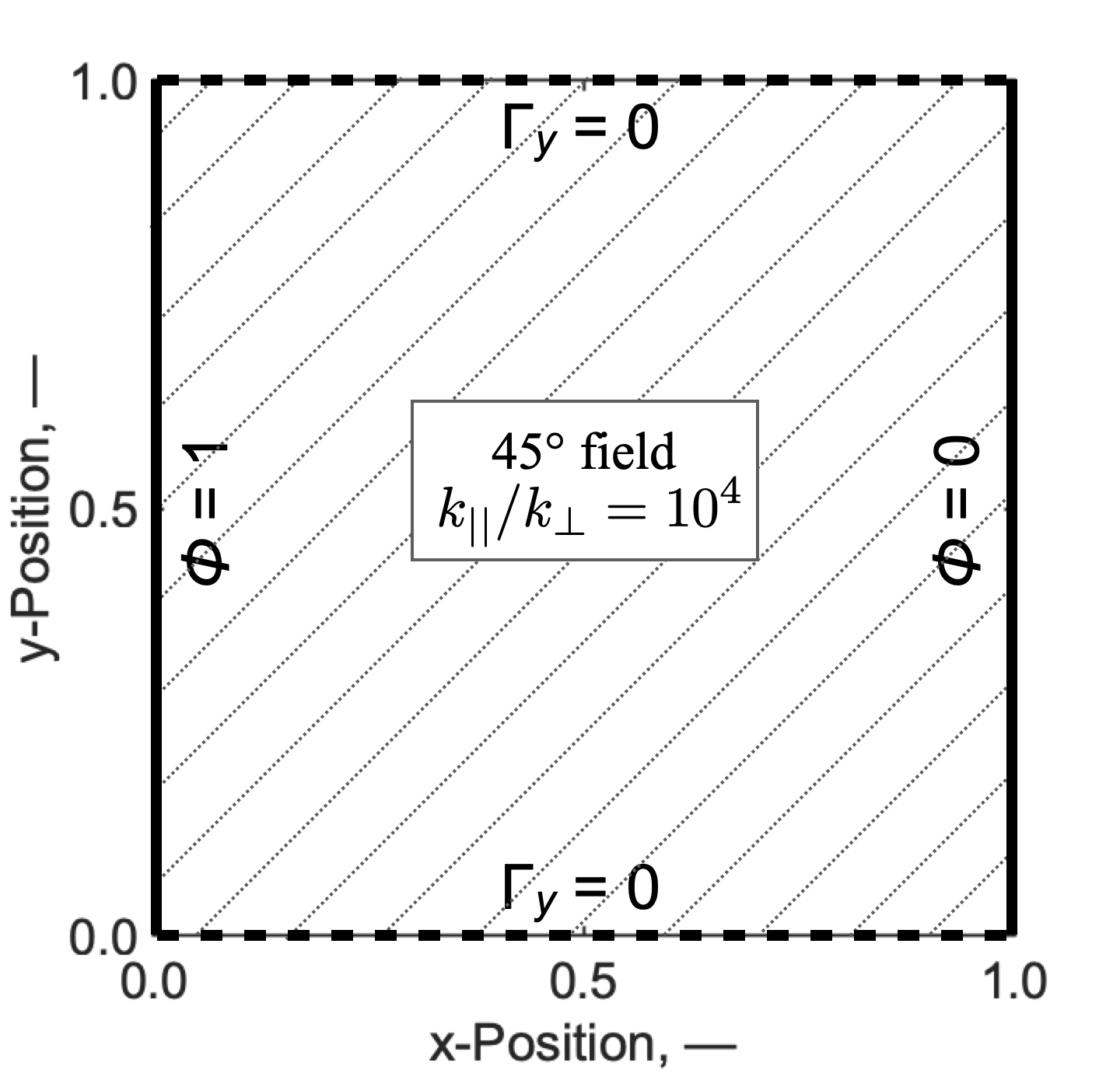}
   		\end{center}
		\vspace{-4mm}
   		\caption{Calculation domain, boundary conditions, and $\pi/4$-inclined magnetic field used for the test.
			Dashed lines denote magnetic field lines with magnetic flux density $B_{\rm ref}$.}
   		\label{fig:condition0}
   \end{figure}

	Another desirable property for the hyperbolic system approach in plasma flow simulation is 3) fulfillment of discrete maximum principle (DMP).
	It is known that standard linear discretizations of anisotropic diffusion with cross-diffusion terms can violate the DMP.
	The DMP is particularly important in plasma flows.
	For example, if the electric potential distribution exhibits spurious peaks in the numerical solution because the DMP is violated, nonphysical behaviors of charged particle motion would be induced \cite{gunter2009}.
	We first show that a linear discretization can indeed fail to satisfy the DMP for anisotropic diffusion.
	\Figure{fig:condition0} illustrates a test case in which a magnetic field inclined at $\pi/4$ is uniformly applied in a $1.0\times1.0$ square domain.
	The left- right-hand boundaries impose Dirichlet conditions for the main variable $\phi$,
	while the top and bottom boundaries impose Neumann conditions with impermeable walls, fixing $\Gamma_y=0$.
	The strength of anisotropy is set to $k_{||}/k_\perp=10^4$.
	This case was computed using a standard diffusion scheme based on second-order central differencing on a 50$\times$50 grid, incorporated with a direct matrix inversion method.
	
\begin{figure}[H]
		\begin{minipage}{0.5\hsize}
			\begin{center}
				\includegraphics[height=55mm]{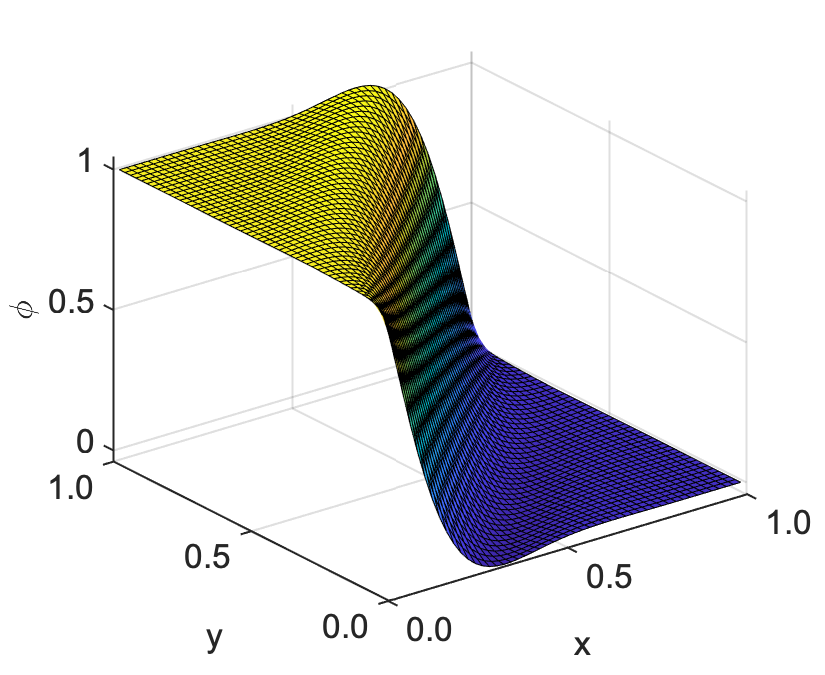}
			\end{center}
			\hspace{40mm}(a)
		\end{minipage}
		\begin{minipage}{0.5\hsize}
			\begin{center}
				\includegraphics[height=55mm]{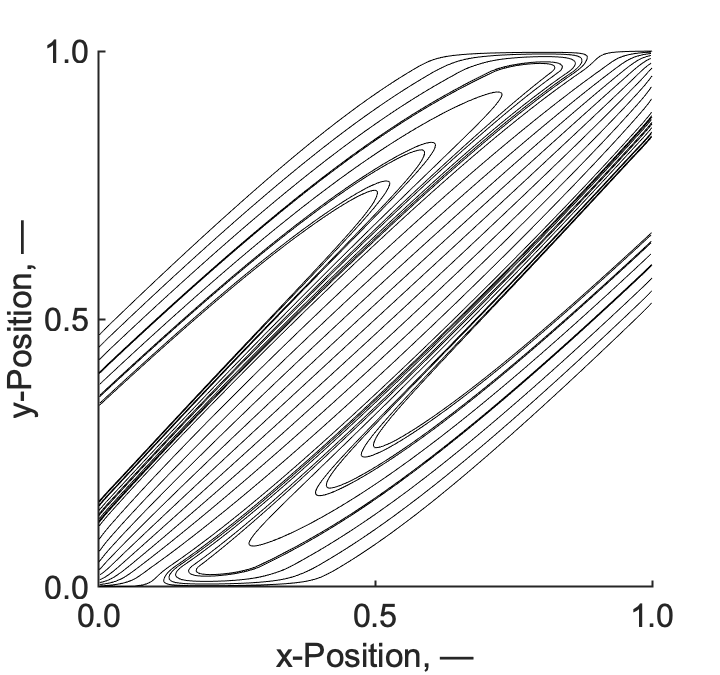}
			\end{center}
			\hspace{40mm}(b)
		\end{minipage}
		\caption{Calculation results of anisotropic diffusion in the test case in \Figure{fig:condition0} . (a) Contours of main variable $\phi$, and (b) streamlines of the flux $\vec{\Gamma}$.}
		\label{fig:failure}
	\end{figure}

	\Figure{fig:failure} shows the numerical results for $\phi$ and the streamlines of $\vec{\Gamma}$.
	As seen in the $\phi$ distribution, the maximum and minimum values exceed the range prescribed by the Dirichlet conditions.
	The streamlines exhibit nonphysical circulations around the peaks of $\phi$.
	Such oscillations can distort simulations in which plasma instabilities or turbulence are involved.
	Similar DMP violations can occur for anisotropic diffusion when the hyperbolic system approach is used.
	Therefore, it is essential to identify conditions under which the hyperbolic system satisfies the DMP.
	
	The issue of DMP violation in anisotropic diffusion and construction of DMP- (or monotonicity-) preserving schemes have been extensively investigated.
	Sharma and Hammett showed that standard algorithms such as central differencing do not preserve monotonicity for anisotropic diffusion and that slope-limiting methods can help mitigate the DMP violations \cite{sharma2007}.
	Several numerical approaches have been proposed to enforce the DMP in anisotropic diffusion, including nonlinear constrained \cite{kuzmin2009}, a posteriori correction \cite{liska2008}, and directional-splitting methods \cite{ngo2016}.
	It is generally understood that nonlinear discretization methods such as a priori or a posteriori limiters are required to satisfy the DMP in diffusion-based (elliptic) schemes.
	The DMP behavior of the hyperbolic system approach must therefore be clarified.
	Chamarthi et al. reported that a fifth-order weighted compact nonlinear scheme (WCNS) combined with a weighted essentially non-oscillatory (WENO) extrapolation is effective in relaxing overshoots and undershoots in results obtained with the hyperbolic system approach \cite{chamarthi2019}.
	However, because lower-order schemes, such as second-order accurate schemes, are widely used in many practical simulations of plasma devices, it remains important to understand the DMP properties even for lower-order discretizations.

	The objective of this study is to assess the DMP properties of the hyperbolic system approach for anisotropic diffusion, and to identify a criterion under which the approach satisfies the DMP.
	We first review the derivation of the hyperbolic system for anisotropic diffusion with a preconditioning chosen to yield a simple system, where the flux Jacobian matrices are identical to the ones in an isotropic diffusion case.
	To obtain a DMP-preserving scheme while maintaining the extensibility and simplicity of the hyperbolic formulation, we avoid nonlinear treatments of numerical fluxes or gradients.
	Instead, we focus on a linear scheme, selecting an appropriate free parameter that appears in the hyperbolic system.
	The role of the parameter is mathematically analyzed in terms of the DMP, and a condition for monotonicity in anisotropic diffusion cases is revealed.
	Fulfillment of the DMP in the hyperbolic system approach with the selected parameter is then verified through several test cases.

\section{Basic equations}\label{sec:basic}
	\subsection{Magnetized electron fluid in quasineutral plasma}
	We consider electron flow in a low-temperature, partially ionized, collisional plasma.
	The following assumptions are commonly made when modeling electron flow in such plasmas:
	\begin{itemize}
		\item Quasineutrality: the plasma is sufficiently dense that the ion and electron number densities are approximately equal (i.e. $n_{\rm e}=n_{\rm i}$).
		\item Inertialess electron: electrons are so light that the inertial effects are negligible owing to frequent collisions with background particles. 
	\end{itemize}
	The first assumption is appropriate for dense plasmas in which the difference between positive-ion and electron densities is small compared with the plasma density \cite{chen2016}.
	The second assumption is valid in plasma sources where a background gas is introduced and a partially ionized gas is generated.
	Under these assumptions, the equation for electrons reduces to a drift-diffusion equation.

	The electron mass conservation equation is
	\begin{equation}
		\nabla\cdot\left(n_{\rm e}\vec{u}_{\rm e}\right)=S_{\rm e},
		\label{eq:mass}
	\end{equation}
	where $n_{\rm e}$ and $u_{\rm e}$ are the electron number density and electron velocity, respectively, and $S_{\rm e}$ is a source term due to chemical reactions such as ionization, recombination, and attachment.
	In quasineutral plasmas, electrons are sufficiently mobile that the electron number density distribution instantaneously becomes equilibrium with the ion number density distribution.
	Hence, the time-derivative terms of electron number density is neglected in the mass conservation equation.
	If one assumes $n_{\rm e} = n_{\rm i}$, Poisson's equation of Gauss's low is not used to obtain the space potential.
	Instead, under the quasineutrality assumption, the space potential is obtained through the electron fluid equations.

	With inertialess electrons, the electron momentum conservation is expressed as
	\begin{equation}
		\nabla\left(en_{\rm e}T_{\rm e}\right) = en_{\rm e}\nabla\phi
		-en_{\rm e} \vec{u}_{\rm e}\times \vec{B}
		-m_{\rm e}\nu_{\rm col}n_{\rm e}\vec{u}_{\rm e},
		\label{eq:momentum1}
	\end{equation}
	where $e$, $T_{\rm e}$, $\phi$, $\vec{B}$, $m_{\rm e}$, and $\nu_{\rm col}$ are the elementary charge, electron temperature, space potential, magnetic field, electron mass, and total collision frequency, respectively.
	From this equation, the electron flux can be written explicitly as
	\begin{equation}
		n_{\rm e}\vec{u}_{\rm e}= n_{\rm e}\left[\mu \right]\nabla \phi
		-\left[\mu\right]\nabla\left(n_{\rm e}T_{\rm e}\right),
		\label{eq:momentum2}
	\end{equation}
	where $\left[\mu\right]$ is the electron mobility tensor.
	The electron flux is driven by gradients of potential and pressure.
	The electron mobility tensor includes the information of magnetic confinement, and it is written in the three dimension as
	\begin{equation}
		\left[\mu\right]=\frac{e}{m_{\rm e}\nu_{\rm col}}\cdot\frac{1}{1+|\vec{\omega}|^2}
			\left[\begin{array}{ccc}
				1+\omega_x^2 & \omega_z+\omega_x\omega_y & -\omega_y+\omega_z\omega_x \\
				-\omega_z+\omega_x\omega_y & 1+\omega_y^2 & \omega_x+\omega_y\omega_z \\
				\omega_y+\omega_z\omega_x & -\omega_x+\omega_y\omega_z & 1+\omega_z^2
			\end{array}\right],
		\label{eq:mobility1}
	\end{equation}
	where $\vec{\omega} = (\omega_x,\omega_y,\omega_z)^T$ is a dimensionless quantity called the Hall parameter, and it is defined for each direction as
	\begin{equation}
		\vec{\omega}=\frac{e}{m_{\rm e}\nu_{\rm col}}\vec{B}.
		\label{eq:omega}
	\end{equation}
	The common definition of the Hall parameter is the ratio of gyrofrequency to collision frequency, and is often directionless.
	However in this paper, the Hall parameter is regarded as directional for convenience of formulation.
	We further assume that the magnetic field is static and that the induced field due to electron current is negligibly small compared with the externally applied magnetic field.
	This assumption is valid for many practical electric propulsion systems and industrial plasma sources.

	In this paper, we consider a two-dimensional problem in the $x-y$ plane where $B_z=0$ with $\vec{B} = (B_x,B_y,B_z)^T$.
	In this case, the electron mobility tensor becomes
	\begin{equation}
		\left[\mu\right]=\frac{e}{m_{\rm e}\nu_{\rm col}}\cdot\frac{1}{1+|\vec{\omega}|^2}
			\left[\begin{array}{cc}
				1+\omega_x^2 & \omega_x\omega_y \\
				\omega_x\omega_y & 1+\omega_y^2 
			\end{array}\right].
		\label{eq:mobility2}
	\end{equation}
	Another useful expression is
	\begin{equation}\label{eq:mobility3}
		\left[\mu\right]=\frac{e}{m_{\rm e}\nu_{\rm col}}\cdot \Theta^{-1}
			\left[\begin{array}{cc}
				1 & \ \\
				\ & \frac{1}{1+|\vec{\omega}|^2}
			\end{array}\right]\Theta,
	\end{equation}
	where $\Theta$ is the rotating matrix defined by the angle between the magnetic field lines and the computational grid.
	This form is conveniently used in practical two-dimensional simulations \cite{komurasaki1995,hagelaar2007}.
		
	The set of equations for the magnetized electron flow consists of Eqs. (\ref{eq:mass}), (\ref{eq:momentum2}), and (\ref{eq:mobility2}).
	To focus on the issues related to anisotropic diffusion, we further assume
	\begin{equation}
		S_{\rm e}=0,\hspace{20pt} n_{\rm e}={\rm const.},\hspace{20pt} T_{\rm e}={\rm const.},\hspace{20pt} \nu_{\rm col}={\rm const.}.
	\end{equation}
	With these assumptions, Eqs. (\ref{eq:mass}), (\ref{eq:momentum2}), and (\ref{eq:mobility2}) reduce to
	\begin{align}
		&\nabla\cdot\vec{\Gamma}=0, \nonumber\\
		&\vec{\Gamma}=-\left[k\right]\nabla\phi,  \nonumber\\
		&\left[k\right]=\left[\begin{array}{cc}
				k_x & k_{\rm c} \\
				k_{\rm c} & k_y
			\end{array}\right]
		=\Theta^{-1}
			\left[\begin{array}{cc}
				1 & \ \\
				\ & \left(k_{||}/k_{\perp}\right)^{-1}
			\end{array}\right]\Theta,
		\label{eq:aniso2}
	\end{align}
	where $\Gamma$ denotes the flux, and $k_{||}/k_{\perp}\equiv 1+|\omega|^2$ determines the magnetic confinement strength; in many magnetically confined plasmas, $k_{||}/k_{\perp} \gg 1$.
	This system captures the essential difficulties of anisotropic diffusion due to magnetic confinement.
	We therefore focus on the problem in \Eqref{eq:aniso2} in what follows.
	
	
	\subsection{Hyperbolic system approach for anisotropic diffusion equation}
	We aim to develop a hyperbolic formulation applicable to anisotropic diffusion problems as well as isotropic diffusion cases.
	We first introduce a hyperbolic system for isotropic diffusion, and then describe the modification of the system for anisotropic diffusion case.
	For isotropic diffusion, the diffusion coefficient is a scaler quantity and the diffusion equation is
	\begin{equation}
		\nabla\cdot\left(-k\nabla \phi\right)=0.
	\end{equation}
	An artificial coefficient $\alpha_{\rm s}$ is introduced as follows:
	\begin{equation}
		\vec{\Gamma}=\alpha_{\rm s}\left(\begin{array}{c}
				u  \\
				v
			\end{array}\right)=-k\nabla \phi.
	\end{equation}
	Here $u$ and $v$ are the variables newly introduced for constructing a hyperbolic system.
	These variables corresponds to the flow velocities. 
	Although one could introduce a vector ($\vec{\alpha}_{\rm s}$) or a tensor ($\left[\vec{\alpha}_{\rm s}\right]$), the parameter is defined as a scaler parameter for simplicity in the present paper.
	
	With the introduced variables, a hyperbolic system can be written as
	\begin{align}
		\begin{split}
			&\frac{\partial \phi}{\partial t} +\frac{\partial}{\partial x}\left(\alpha_{\rm s}u\right)+\frac{\partial}{\partial y}\left(\alpha_{\rm s}v\right)=0,\\
			&\frac{\partial u}{\partial t}+k\frac{\partial \phi}{\partial x}=-\alpha_{\rm s}u, \\
			&\frac{\partial v}{\partial t}+k\frac{\partial \phi}{\partial y}=-\alpha_{\rm s}v,
		\end{split}
	\end{align}
	where $t$ is a pseudo time introduced for iterative calculation.
	This equation system is equivalent to the original diffusion equation once the pseudo-time derivative terms become negligibly small at steady state.
	Further, one can consider multiplying a preconditioning matrix $P$ to the pseudo-time derivative terms.
	For the isotropic case, the preconditioning matrix is chosen as
	\begin{equation}
		P=\left[\begin{array}{ccc}
			\alpha_{\rm s}  & \ & \ \\
			\ & k  & \ \\
			\ & \ &k \\
		\end{array}\right],
	\end{equation}
	which yields	
	\begin{align}\label{eq:set_iso}
		&\frac{\partial \phi}{\partial t} +\frac{\partial u}{\partial x}+\frac{\partial v}{\partial y}=0, \nonumber\\
		&\frac{\partial u}{\partial t}+\frac{\partial \phi}{\partial x}\hspace{23pt}=-\frac{\alpha_{\rm s}}{k}u, \\
		&\frac{\partial v}{\partial t}\hspace{23pt}+\frac{\partial \phi}{\partial y}=-\frac{\alpha_{\rm s}}{k}v. \nonumber
	\end{align}
	With this system, the Jacobian matrices for the numerical fluxes in the x- and y-directions are written as
	\begin{equation}
		A=
		\left[\begin{array}{ccc}
			0 & 1 & 0 \\
			1 & 0 & 0 \\
			0 & 0 & 0 \\
		\end{array}\right],\hspace{15pt}
		B=
		\left[\begin{array}{ccc}
			0 & 0 & 1 \\
			0 & 0 & 0 \\
			1 & 0 & 0 \\
		\end{array}\right].
		\label{eq:Jacobi}
	\end{equation}
	The eigenvalues of these matrices are $\lambda_A=-1,0,1$ and $\lambda_B=-1,0,1$.
	Based on these Jacobian matrices, one can construct an upwind method based on the approximate Riemann solver.
	
	For the anisotropic diffusion case involving the diffusion tensor in \Eqref{eq:aniso2}, the gradient variables can be defined as
	\begin{equation}
		\alpha_{\rm s}\left(\begin{array}{c}
				u  \\
				v
			\end{array}\right)=-\left[\begin{array}{cc}
				k_x & k_{\rm c} \\
				k_{\rm c} & k_y
			\end{array}\right]\nabla \phi,
	\end{equation}
	which yields the hyperbolic system
	\begin{align}\label{eq:hyperbolic}
		&\frac{\partial \phi}{\partial t} +\alpha_{\rm s}\frac{\partial u}{\partial x}+\alpha_{\rm s}\frac{\partial v}{\partial y}=0, \nonumber\\
		&\frac{\partial u}{\partial t}+k_x\frac{\partial \phi}{\partial x}+k_c\frac{\partial \phi}{\partial y}
		=-\alpha_{\rm s}u,\nonumber \\
		&\frac{\partial v}{\partial t}+k_c\frac{\partial \phi}{\partial x}+k_y\frac{\partial \phi}{\partial y}
		=-\alpha_{\rm s}v. \nonumber\\
	\end{align}
	A preconditioning matrix is also introduced for this hyperbolic system.
	In this paper, we use a preconditioner as
	\begin{equation}
		P=\left[\begin{array}{ccc}
			\alpha_{\rm s}  & \ & \ \\
			\ & k_x  & k_{\rm c} \\
			\ & k_{\rm c} &k_y \\
		\end{array}\right].
	\end{equation}
	Then, the hyperbolic system is modified as
	\begin{align}\label{eq:set_aniso}
		\begin{split}
			&\frac{\partial \phi}{\partial t} +\frac{\partial u}{\partial x}+\frac{\partial v}{\partial y}=0,\\
			&\frac{\partial u}{\partial t}+\frac{\partial \phi}{\partial x}
				\hspace{23pt}=-\frac{k_y\alpha_{\rm s}}{\Delta}u+\frac{k_c\alpha_{\rm s}}{\Delta}v, \\
			&\frac{\partial v}{\partial t}
			\hspace{23pt}+\frac{\partial \phi}{\partial y}
			= \frac{k_c\alpha_{\rm s}}{\Delta}u-\frac{k_x\alpha_{\rm s}}{\Delta}v,
		\end{split}
	\end{align}
	where $\Delta \equiv k_xk_y-k_c^2$, which is positive in the cases of the anisotropy due to magnetic confinement.
	Based on the tensor in \Eqref{eq:aniso2}, $\Delta$ can be evaluated as
	\begin{equation}
		\Delta=\frac{1}{\left(1+|\vec{\omega}|^2\right)^2}\left(1+\omega_x^2+\omega_y^2\right)>0.
	\end{equation}
	
	The effects of anisotropy are included in the source terms on the right-hand side.
	Defining the source term Jacobian
	\begin{equation}
		H=\left[\begin{array}{ccc}
			\  & \ & \ \\
			\ & \frac{k_y\alpha_{\rm s}}{\Delta} & -\frac{k_{\rm c}\alpha_{\rm s}}{\Delta} \\
			\ &-\frac{k_{\rm c}\alpha_{\rm s}}{\Delta} & -\frac{k_{\rm x}\alpha_{\rm s}}{\Delta} \\
		\end{array}\right],
	\end{equation}
	the source term can be written as $S=-HU$.
	With the variable vector $U\equiv \left(\phi,u,v\right)^T$, the system can be expressed as
	\begin{equation}\label{eq:hyperbolic2}
		\frac{\partial U}{\partial t}+A\frac{\partial U}{\partial x}+B\frac{\partial U}{\partial y}=-HU.
	\end{equation}
	The flux Jacobian matrices $A$ and $B$ are identical to those in the isotropic case in \Eqref{eq:Jacobi}.
	Therefore, the numerical fluxes associated with the space derivatives can be treated exactly as in the isotropic diffusion case.
	One just needs to modify the source term treatment to extend the hyperbolic system approach from the isotropic to anisotropic diffusion case.
	
\section{Monotonicity of FDM}\label{sec:mon}
\subsection{Viscosity solutions}
In this section, we will explain that the monotonicity of the finite difference method (FDM)
for a class of partial differential equations is important for its convergence property.
This class of partial differential equations admits a notion of weak solutions called \textit{viscosity solutions}.
We note that solving \Eqref{eq:ade} is equivalent to find a stationary solution of the following system of partial differential equations:
\begin{equation}\label{eq:ade_st}
	\begin{cases}
		\phi_t - \nabla\cdot(\matK\nabla\phi) = 0&\qquad\mbox{in} \quad\Omega_T := \Omega\times(0,T),\\
		\phi(\cdot,0) = \phi_0&\qquad\mbox{in} \quad\Omega,
	\end{cases}
\end{equation}
for some given initial data $\phi_0$.
Here, we briefly recall the notion of viscosity solutions for \Eqref{eq:ade_st}.
The first equation of the above system can be rewritten as
\begin{equation}\label{eq:gf}
    \phi_t + F(x,t,\phi(x,t),\nabla\phi(x,t),\nabla^2\phi(x,t)) = 0\qquad\mbox{in}\quad \Omega\times(0,T),
\end{equation}
where $\nabla^2\phi$ denotes the Hessian matrix of $\phi$, say $\nabla^2\phi = \left[\partial_{ij}\phi\right]_{i,j}$ for $1\leq i,\,j\leq 2$, and
$F:\closure{\Omega}\times[0,T]\times\mathbb{R}\times\mathbb{R}^2\times\mathbb{S}^2\to\mathbb{R}$ with $\mathbb{S}^2$ being the set of all $2\times 2$ symmetric matrices.
To represent the anisotropic diffusion equation \eqref{eq:ade_st} in the form of \eqref{eq:gf}, we set
\begin{equation}\label{eq:gf_ade}
    F(x,t,r,\vec{p},X) := -\matK:\mat{X},
\end{equation}
where $\mat{X}:\mat{Y}$ denotes the Frobenius inner product of matrices $\mat{X}$ and $\mat{Y}$, say $\mat{X}:\mat{Y} = \sum_{i,j}x_{ij}y_{ij}$ for $\mat{X} = \mat{x_{ij}}\in\mathbb{R}^{2\times 2}$ and $\mat{Y} = \mat{y_{ij}}\in\mathbb{R}^{2\times 2}$.
This inner product can be rewritten as $\mat{X}:\mat{Y} = \trace{\mat{X}^T\mat{Y}}$, where $\trace{\mat{X}}$ denotes the trace of $\mat{X}$, and $\mat{X}^T$ is the transpose of $\mat{X}$.

To establish the well-posedness of \Eqref{eq:gf} with regard to viscosity solutions, we suppose that $F$ is \textit{degenerate elliptic}.
Precisely speaking, we assume that

\begin{description}
	\item[(DE)] for every $\mat{X}$, $\mat{Y}\in\mathbb{S}^2$ with $\mat{X}\geq \mat{Y}$, we have
	\begin{equation*}
		F(x,t,r,\vec{p},\mat{X}) \leq F(x,t,r,\vec{p},\mat{Y})\quad\mbox{for all}\quad (x,t,r,\vec{p})\in\closure{\Omega}\times[0,T]\times\mathbb{R}\times\mathbb{R}^2,
	\end{equation*}
	where $\mat{X}\geq \mat{Y}$ means that $\mat{X} -\mat{Y}$ is a positive semi-definite matrix.
\end{description}

Then, we define viscosity solutions for \eqref{eq:gf} as follows:

\begin{dfn}[Viscosity solutions]
    A function $u$ is called a viscosity sub- (resp. super-) solution to \eqref{eq:gf}
    provided that for every smooth function $\varphi^+:\Omega_T\to\mathbb{R}$ (resp. $\varphi^-:\Omega_T\to\mathbb{R}$)
    such that $u - \varphi^+$ (resp. $u - \varphi^-$) has a local maximum (resp. a local minimum) at $(x,t)\in\Omega_T$, it follows that
    \begin{align*}
        \varphi^+_t + F(x,t,\varphi^+(x,t),\nabla\varphi^+(x,t),\nabla^2\varphi^+(x,t)) &\leq 0.\\
        (resp.\ \varphi^-_t + F(x,t,\varphi^-(x,t),\nabla\varphi^-(x,t),\nabla^2\varphi^-(x,t)) &\geq 0.)
    \end{align*}
\end{dfn}

\begin{rem}
It is easily seen that a classical solution to \eqref{eq:gf} is a viscosity solution when $F$ is degenerate elliptic.
\end{rem}

We now recall a comparison principle of viscosity solutions for the parabolic equation \eqref{eq:gf}:

\begin{prop}[Theorem 3.1.1 of \cite{G}]\label{prop:cp}
    Assume that $\Omega$ is bounded, and $F$ is degenerate elliptic. Assume further that $F$ satisfies the following conditions:
    \begin{description}
        \item[(F1)] $F$ is continuous in $\closure{\Omega}\times[0,T]\times\mathbb{R}\times\mathbb{R}^2\times\mathbb{S}^2$.
        \item[(F2)] There exists a constant $c_0\in \mathbb{R}$ such that the mapping
        \begin{equation*}
            \mathbb{R}\ni r\mapsto F(x,t,r,\vec{p},\mat{X}) + c_0r\in\mathbb{R}
        \end{equation*}
        is a non-decreasing function for every $(x,t,\vec{p},\mat{X})\in\closure{\Omega}\times[0,T]\times\mathbb{R}^2\times\mathbb{S}^2$.
    \end{description}
    Assume that $u$ and $v$ are viscosity sub- and supersolution to \eqref{eq:gf} in $\Omega_T$.
    Then, it holds that $u\leq v$ in $\Omega_T$ provided that $u\leq v$ on $\partial_p\Omega_T := (\pOmega\times(0,T))\cup(\Omega\times\{0\})$.
    In particular, if $u$ and $v$ are viscosity solutions to \eqref{eq:gf} with the initial data $u_0$ and $v_0$ and the Dirichlet data $u_D$ and $v_D$, respectively.
    Then, $u\leq v$ in $Q$ holds provided that $u_0\leq v_0$ in $\Omega$ and $u_D\leq v_D$ on $\pOmega\times(0,T)$.
\end{prop}
\begin{rem}
    The original statement in \cite[Theorem 3.1.1]{G} applies to more general $F$.
    Precisely speaking, the function $F$ is allowed not to be defined at $\vec{p} = \vec{0}$, although
    we need to assume that $F_*(x,t,r,\vec{0},O) = F^*(x,t,r,\vec{0},O)$ is finite,
    where $F_*$ and $F^*$ are the lower and upper semicontinuous envelopes of $F$, respectively.
    In other words, if $F$ can be extended to a continuous function on $\closure{\Omega}\times[0,T]\times\mathbb{R}\times\mathbb{R}^2\times\mathbb{S}^2$,
    then the statement of \Proposition{prop:cp} still holds.
\end{rem}
For the choice \eqref{eq:gf_ade} of $F$ and $c_0 \geq 0$,
it is easily observed that $F$ is degenerate elliptic, and the conditions (F1) and (F2) are satisfied since $\matK$ is a positive definite matrix,
and hence we can apply \Proposition{prop:cp} to the anisotropic diffusion equation \eqref{eq:ade_st}.
The comparison principle for viscosity solutions of anisotropic diffusion equations states
that the order of the initial functions $u_0\leq v_0$ is preserved until the time horizon $T$.
This fact motivates us to establish the monotonicity of the finite difference scheme for \eqref{eq:ade_st}.

We now recall sufficient conditions to ensure the convergence of a discrete scheme $\discreteScheme{\Delta t}{}$
from the seminal work by Barles and Souganidis \cite{BarlesSouganidis1991}:

\begin{prop}[Theorem 2.1 of \cite{BarlesSouganidis1991}]\label{prop:bs}
    Let $UC(\closure{\Omega})$ be the set of all uniformly continuous functions on $\closure{\Omega}$.
    Assume that a family of function operators $\{\discreteScheme{\Delta t}{}\}_{\Delta t > 0}$
    between $UC(\closure{\Omega})$ and itself satisfies the following conditions:\\\\
    \noindent
    \begin{subequations}\label{eq:bs}
    \textbf{[Monotonicity]} For every $u,\,v\in UC(\closure{\Omega})$ with $u\leq v$, it follows that\\
    \begin{equation}\label{eq:bs_mon}
        \discreteScheme{\Delta t}{}u \leq \discreteScheme{\Delta t}{} v.
    \end{equation}
    \textbf{[Translation invariance]} For every $u\in UC(\closure{\Omega})$ and $c\in\mathbb{R}$, it follows that\\
    \begin{equation}\label{eq:bs_ti}
        \discreteScheme{\Delta t}{}(u+c) = \discreteScheme{\Delta t}{}u + c.
    \end{equation}
    \textbf{[Consistency]} For every smooth function $\varphi$ in $\closure{\Omega}$, it follows that\\
    \begin{equation}\label{eq:bs_con}
        \lim_{\Delta t\to 0}\frac{\discreteScheme{\Delta t}{}\varphi - \varphi}{\Delta t} = \matK:\nabla^2\varphi.
    \end{equation}
    \end{subequations}
    Let an approximate solution $u^{\Delta t}(x,t)$ be defined by
    \begin{equation*}
        u^{\Delta t}(x,t) := \discreteScheme{\Delta t}{\lfloor\frac{t}{\Delta t}\rfloor}u_0(x)\quad\mbox{for all}\quad (x,t)\in\closure{\Omega}\times[0,T],
    \end{equation*}
    where $u_0(x)$ is given initial function, and $\lfloor s\rfloor$ denotes the largest integer less than or equal to $s$ for $s\in\mathbb{R}$.
    Then, $u^{\Delta t}$ uniformly converges to the viscosity solution $u$ of \eqref{eq:gf} in $\closure{\Omega}\times(0,T)$ as $\Delta t\to 0$
    whenever the target equation \eqref{eq:gf} satisfies the comparison principle.
\end{prop}
We see that typical FDMs satisfy the translation invariance. Moreover, we can expect that the consistency condition holds for FDMs approximating \eqref{eq:hyperbolic} since $\alpha_s$ is constant.
Therefore, we focus on the monotonicity of FDM which is one of the key ingredients to ensure the convergence of the numerical method as stated in \Proposition{prop:bs}.
We note that FDM approximating viscosity solutions to anisotropic diffusion equation \eqref{eq:ade_st} preferably possesses a comparison principle according to \Proposition{prop:cp}.
In fact, the comparison principle for viscosity solutions of anisotropic diffusion equations also implies a kind of DMP for the FDM approximating \eqref{eq:ade_st}.

\begin{cor}\label{cor:dmp}
	Assume that a discrete scheme $\discreteScheme{\Delta t}{}$ satisfies the conditions \eqref{eq:bs_mon}, \eqref{eq:bs_ti}, and \eqref{eq:bs_con} in \Proposition{prop:bs}.
	Then, for every $\varepsilon > 0$, there exists a $\Delta t_0 > 0$ such that for every $\Delta t < \Delta t_0$, the following holds:
	\begin{equation}\label{eq:dmp}
		m - \varepsilon \leq \discreteScheme{\Delta t}{\lfloor\frac{t}{\Delta t}\rfloor}u_0(x) \leq M + \varepsilon\qquad\mbox{for all}\quad x\in\closure{\Omega}\quad \mbox{and} \quad t\in(0,T),
	\end{equation}
	where $m := \min_{\closure{\Omega}}{u_D}$ and $M := \max_{\closure{\Omega}}{u_D}$.
\end{cor}
\begin{proof}
	We deduce from \Proposition{prop:bs} that $u^{\Delta t}(x,t) = \discreteScheme{\Delta t}{\lfloor\frac{t}{\Delta t}\rfloor}u_0(x)$ uniformly converges to the unique viscosity solution $u(x,t)$ of \eqref{eq:ade_st} as $\Delta t\to 0$.
	Meanwhile, we see that $u\equiv m$ and $u\equiv M$ are classical sub- and supersolutions to \eqref{eq:ade_st} with the initial data $u_0$ and the Dirichlet data $u_D$.
	Thus, the comparison principle for viscosity solutions (\Proposition{prop:cp}) implies that
	$m\leq u(x,t)\leq M$ for all $(x,t)\in\closure{\Omega}\times(0,T)$. Therefore, \Eqref{eq:dmp} follows for sufficiently small $\Delta t > 0$.
\end{proof}

According to \Corollary{cor:dmp}, if $t$ is pseudo-time, then we can take $\Delta t > 0$ so small that
the discrete solution $\discreteScheme{\Delta t}{\lfloor\frac{t}{\Delta t}\rfloor}u_0(x)$ is nearly bounded by the Dirichlet data $u_D$.
However, as discussed in \Section{sec:intro}, the hyperbolic system approach for anisotropic diffusion problems can cause a failure of DMP (see \Figure{fig:failure}),
and we are led to consider a monotonicity-preserving scheme for the hyperbolic system \Eqref{eq:hyperbolic}.
\subsection{Monotonicity of the hyperbolic system}
As a beginning point, we observe that a naive FDM for \eqref{eq:ade_st}
is not monotone with respect to the order of the initial data.
Let $\Omega$ be the square domain $(0,1)^2\subset\mathbb{R}^2$.
For $N_x,\,N_y \geq 2$, we set $x_{i,j} := (i\Delta x, j\Delta y)$ for $0\leq i\leq N_x$ and $0\leq j\leq N_y$ with $\Delta x := 1/N_x$ and $\Delta y := 1/N_y$.
We now approximate a solution $\phi$ to \eqref{eq:ade_st} by the central difference approximation:
\begin{align}\label{eq:ade_cls}
    \begin{split}
    \discreteFunc{\phi}{i}{j}{n+1} &=\left(1 - 2k_x\frac{\Delta t}{\Delta x^2} - 2k_y\frac{\Delta t}{\Delta y^2}\right)\discreteFunc{\phi}{i}{j}{n} + k_x\frac{\Delta t}{\Delta x^2}\left(\discreteFunc{\phi}{i+1}{j}{n} + \discreteFunc{\phi}{i-1}{j}{n}\right)\\
    &+ k_y\frac{\Delta t}{\Delta y^2}\left(\discreteFunc{\phi}{i}{j+1}{n} + \discreteFunc{\phi}{i}{j-1}{n}\right)
    + \frac{k_c}{2}\frac{\Delta t}{\Delta x\Delta y}\left(\discreteFunc{\phi}{i+1}{j+1}{n} - \discreteFunc{\phi}{i+1}{j-1}{n} - \discreteFunc{\phi}{i-1}{j+1}{n} + \discreteFunc{\phi}{i-1}{j-1}{n}\right),
    \end{split}
\end{align}
where $\discreteFunc{\phi}{i}{j}{n}$ denotes a discretization of $\phi(x_{i,j},n\Delta t)$ for each time step $n \geq 0$.
The formula \eqref{eq:ade_cls} defines a function operator $\discreteScheme{\Delta t}{}$,
although the scheme also depends on the choice of the spatial mesh sizes $\Delta x$ and $\Delta y$.
We easily see that $\discreteScheme{\Delta t}{}$ is monotone in the meaning of \Eqref{eq:bs_mon} if all coefficients of $\discreteFunc{\varphi}{\ell}{k}{b}$'s are all non-negative.
However, due to the cross terms, this $\discreteScheme{\Delta t}{}$ is not monotone with respect to the order of the previous data.
Given the previous data $\discreteFunc{\phi}{i}{j}{n}$, $\discreteFunc{u}{i}{j}{n}$ and $\discreteFunc{v}{i}{j}{n}$,
we discretize the system \eqref{eq:hyperbolic2} with a time step $\Delta t$ and mesh sizes $\Delta x$ and $\Delta y$, and we obtain
\begin{equation*}
	\frac{\discFunc{U}{n+1} - \discFunc{U}{n}}{\Delta t}
	= -\frac{A^+}{\Delta x}\left(\discreteFunc{U}{i}{j}{n} - \discreteFunc{U}{i-1}{j}{n}\right)
	-\frac{A^-}{\Delta x}\left(\discreteFunc{U}{i+1}{j}{n} - \discreteFunc{U}{i}{j}{n}\right)
	-\frac{B^+}{\Delta y}\left(\discreteFunc{U}{i}{j}{n} - \discreteFunc{U}{i}{j-1}{n}\right)
	-\frac{B^-}{\Delta y}\left(\discreteFunc{U}{i}{j+1}{n} - \discreteFunc{U}{i}{j}{n}\right)
	-H\discreteFunc{U}{i}{j}{n+1},
\end{equation*}
where
\begin{equation*}
	A^+ := \left[
		\begin{array}{ccc}
			\frac{1}{2} & \frac{1}{2} & \\
			\frac{1}{2} & \frac{1}{2} & \\
			 &  & 
		\end{array}
	\right],\quad
	A^- := \left[
		\begin{array}{ccc}
			-\frac{1}{2} & \frac{1}{2} & \\
			\frac{1}{2} & -\frac{1}{2} & \\
			 &  & 
		\end{array}	
		\right],\quad
	B^+ := \left[
		\begin{array}{ccc}
			\frac{1}{2} & & \frac{1}{2} \\
			 &  & \\
			\frac{1}{2} & & \frac{1}{2} \\
		\end{array}
		\right],\quad\mbox{and}\quad
	B^- := \left[
		\begin{array}{ccc}
			-\frac{1}{2} & & \frac{1}{2} \\
			 &  & \\
			\frac{1}{2} & & -\frac{1}{2} \\
		\end{array}
		\right].
\end{equation*}
Solving the above equation with respect to $\discreteFunc{U}{i}{j}{n+1}$, we have the following FDM:
\begin{align}\label{eq:ade_st3}
    \begin{split}
    &\discreteFunc{\phi}{i}{j}{n+1} = \left(1-\frac{\Delta t}{\Delta x} - \frac{\Delta t}{\Delta y}\right)\discreteFunc{\phi}{i}{j}{n} + \frac{\Delta t}{2\Delta x}\discreteFunc{\phi}{i-1}{j}{n} + \frac{\Delta t}{2\Delta x}\discreteFunc{\phi}{i+1}{j}{n} + \frac{\Delta t}{2\Delta y}\discreteFunc{\phi}{i}{j-1}{n} + \frac{\Delta t}{2\Delta y}\discreteFunc{\phi}{i}{j+1}{n}\\
    &\quad +\frac{\Delta t}{2\Delta x}\discreteFunc{u}{i-1}{j}{n} - \frac{\Delta t}{2\Delta x}\discreteFunc{u}{i+1}{j}{n} + \frac{\Delta t}{2\Delta y}\discreteFunc{v}{i}{j-1}{n} - \frac{\Delta t}{2\Delta y}\discreteFunc{v}{i}{j+1}{n},
    \end{split}
\end{align}
\begin{align*}
    (1 + \beta_y)\discreteFunc{u}{i}{j}{n+1} - \beta_c\discreteFunc{v}{i}{j}{n+1} &= \left(1-\frac{\Delta t}{\Delta x}\right)\discreteFunc{u}{i}{j}{n} + \frac{\Delta t}{2\Delta x}\discreteFunc{u}{i-1}{j}{n} + \frac{\Delta t}{2\Delta x}\discreteFunc{u}{i+1}{j}{n} + \frac{\Delta t}{2\Delta x}\discreteFunc{\phi}{i-1}{j}{n}- \frac{\Delta t}{2\Delta x}\discreteFunc{\phi}{i+1}{j}{n},\\
    -\beta_c\discreteFunc{u}{i}{j}{n+1} + \left(1 + \beta_x\right)\discreteFunc{v}{i}{j}{n+1} &= \left(1 - \frac{\Delta t}{\Delta y}\right)\discreteFunc{v}{i}{j}{n} + \frac{\Delta t}{2\Delta y}\discreteFunc{v}{i}{j-1}{n} + \frac{\Delta t}{2\Delta y}\discreteFunc{v}{i}{j+1}{n} + \frac{\Delta t}{2\Delta y}\discreteFunc{\phi}{i}{j-1}{n} - \frac{\Delta t}{2\Delta y}\discreteFunc{\phi}{i}{j+1}{n},
\end{align*}
where
\begin{equation*}
    \beta_x := \frac{k_x\alpha_s\Delta t}{\detK},\quad \beta_y := \frac{k_y\alpha_s\Delta t}{\detK},\quad\mbox{and}\quad \beta_c := \frac{k_c\alpha_s\Delta t}{\detK}.
\end{equation*}
Solving the last two equalities of the above system with respect to $\discreteFunc{u}{i}{j}{n+1}$ and $\discreteFunc{v}{i}{j}{n+1}$, we have
\begin{equation}\label{eq:ade_st4}
	\left[
		\begin{array}{c}
        	\discreteFunc{u}{i}{j}{n+1}\\
        	\discreteFunc{v}{i}{j}{n+1}
		\end{array}
	\right] = 
	\left[
		\begin{array}{cc}
        	\frac{1+\beta_x}{|B|} & \frac{\beta_c}{|B|}\\
        	\frac{\beta_c}{|B|} & \frac{1+\beta_y}{|B|}
		\end{array}
	\right]
	\left[
		\begin{array}{c}
        \left(1-\frac{\Delta t}{\Delta x}\right)\discreteFunc{u}{i}{j}{n} + \frac{\Delta t}{2\Delta x}\discreteFunc{u}{i-1}{j}{n} + \frac{\Delta t}{2\Delta x}\discreteFunc{u}{i+1}{j}{n} + \frac{\Delta t}{2\Delta x}\discreteFunc{\phi}{i-1}{j}{n} - \frac{\Delta t}{2\Delta x}\discreteFunc{\phi}{i+1}{j}{n}\\
        \left(1-\frac{\Delta t}{\Delta y}\right)\discreteFunc{v}{i}{j}{n} + \frac{\Delta t}{2\Delta y}\discreteFunc{v}{i}{j-1}{n} + \frac{\Delta t}{2\Delta y}\discreteFunc{v}{i}{j+1}{n} + \frac{\Delta t}{2\Delta y}\discreteFunc{\phi}{i}{j-1}{n} - \frac{\Delta t}{2\Delta y}\discreteFunc{\phi}{i}{j+1}{n}
		\end{array}
	\right],
\end{equation}
where
\begin{equation*}
    B := \left[
		\begin{array}{cc}
        	1+\beta_y & -\beta_c\\
        	-\beta_c & 1+\beta_x
		\end{array}
	\right].
\end{equation*}
Here, we note that $|B|$ is positive (resp. negative) if and only if either $\alpha_s < \alpha_s^-$ or $\alpha_s^+ < \alpha_s$ (resp. $\thresM < \alpha_s < \thresP$), where
\begin{equation}\label{eq:thres}
    \alpha_s^\pm := \frac{-(k_x + k_y)\pm\sqrt{(k_x-k_y)^2 + 4k_c^2}}{2\Delta t}.
\end{equation}
Thus, the system \eqref{eq:ade_st4} is solvable,
that is $B$ is invertible, if neither $\alpha_s = \alpha_s^-$ nor $\alpha_s = \alpha_s^+$.
Seeing the first step of the scheme \eqref{eq:ade_st3}, it seems obviously monotone with respect to the previous step of $\phi$.
However, the next step $\discreteFunc{\phi}{i}{j}{n+1}$ of $\discreteFunc{\phi}{i}{j}{n}$
suffers from the previous steps of $\discreteFunc{\phi}{i}{j}{k}$ for $k \leq n-1$
regarding the monotonicity of the scheme. To observe the impact of the parameter $\alpha_s$,
let us translate the scheme \eqref{eq:ade_st3} into the one of $\phi$.
We use the formula \eqref{eq:ade_st4} to compute the values of $\discreteFunc{u}{i-1}{j}{n}$, $\discreteFunc{u}{i+1}{j}{n}$, $\discreteFunc{v}{i}{j-1}{n}$ and $\discreteFunc{v}{i}{j+1}{n}$
on replacing $n$ by $n-1$. For simplicity, we suppose that $\discreteFunc{u}{i}{j}{n-1} = \discreteFunc{v}{i}{j}{n-1}\equiv 0$.
Then, it holds that
\begin{equation}\label{eq:ade_st6}
    \discreteFunc{\phi}{i}{j}{n+1} = \left(1-\frac{\Delta t}{\Delta x} - \frac{\Delta t}{\Delta y}\right)\discreteFunc{\phi}{i}{j}{n} + \frac{\Delta t}{2\Delta x}\discreteFunc{\phi}{i-1}{j}{n} + \frac{\Delta t}{\Delta x}\discreteFunc{\phi}{i+1}{j}{n} + \frac{\Delta t}{2\Delta y}\discreteFunc{\phi}{i}{j-1}{n} + \frac{\Delta t}{2\Delta y}\discreteFunc{\phi}{i}{j+1}{n}
    + \sum_{\substack{i-2\leq k\leq i+2\\ j-2\leq \ell\leq j+2}}\discreteFunc{C}{k}{\ell}{n-1}\discreteFunc{\phi}{k}{\ell}{n-1},
\end{equation}
where 
\begin{align*}
    \discreteFunc{C}{i-1}{j-1}{n-1} = \frac{2\beta_c}{|B|}\left(\frac{\Delta t}{2\Delta x}\right)\left(\frac{\Delta t}{2\Delta y}\right),
    & \qquad\discreteFunc{C}{i+1}{j+1}{n-1} = -\frac{2\beta_c}{|B|}\left(\frac{\Delta t}{2\Delta x}\right)\left(\frac{\Delta t}{2\Delta y}\right),\\
    \discreteFunc{C}{i-2}{j}{n-1} = \left(\frac{\Delta t}{2\Delta x}\right)^2\frac{1+\beta_x}{|B|},
    & \qquad\discreteFunc{C}{i}{j-2}{n-1} =  \left(\frac{\Delta t}{2\Delta y}\right)^2\frac{1+\beta_y}{|B|},\\
    \discreteFunc{C}{i+2}{j}{n-1} = -\left(\frac{\Delta t}{2\Delta x}\right)^2\frac{1+\beta_x}{|B|},
    & \qquad\discreteFunc{C}{i}{j+2}{n-1} =  -\left(\frac{\Delta t}{2\Delta y}\right)^2\frac{1+\beta_y}{|B|}
\end{align*}
and $\discreteFunc{C}{k}{\ell}{n-1} = 0$ for the other pairs $(k,\ell)$.
Due to these cross terms of the $(n-1)$'s step,
the scheme \eqref{eq:ade_st3} is not monotone with respect to the previous steps of $\phi$ unless $\beta_c = 0$, which is equivalent to either $k_c=0$ or $\alpha_s=0$;
both cases are not of interest in the context of anisotropic diffusion equations.
As stated in \Proposition{prop:bs},
the monotonicity of discrete schemes is crucial to ensure the convergence of the numerical method.
Therefore, we shall modify FDM \eqref{eq:ade_st3} so that its monotonicity property is nearly preserved.
Explicitly, we introduce the following scheme:
\begin{multline}\label{eq:ade_st5}
    \discreteFunc{\phi}{i}{j}{n+1} := \left(1-\frac{\Delta t}{\Delta x} - \frac{\Delta t}{\Delta y}\right)\discreteFunc{\phi}{i}{j}{n} + \frac{\Delta t}{2\Delta x}\discreteFunc{\phi}{i-1}{j}{n} + \frac{\Delta t}{2\Delta x}\discreteFunc{\phi}{i+1}{j}{n} + \frac{\Delta t}{2\Delta y}\discreteFunc{\phi}{i}{j-1}{n} + \frac{\Delta t}{2\Delta y}\discreteFunc{\phi}{i}{j+1}{n}\\
    + \frac{\Delta t}{2\Delta x}\discreteFunc{u}{i-1}{j}{n+1} - \frac{\Delta t}{2\Delta x}\discreteFunc{u}{i+1}{j}{n+1} + \frac{\Delta t}{2\Delta y}\discreteFunc{v}{i}{j-1}{n+1} - \frac{\Delta t}{2\Delta y}\discreteFunc{v}{i}{j+1}{n+1}
\end{multline}
with the assumption that $\discreteFunc{u}{i}{j}{n} = \discreteFunc{v}{i}{j}{n}\equiv 0$.
In contrast to the scheme \eqref{eq:ade_st3}, the next steps $\discreteFunc{u}{i}{j}{n+1}$ and $\discreteFunc{v}{i}{j}{n+1}$ of $\discreteFunc{u}{i}{j}{n}$ and $\discreteFunc{v}{i}{j}{n}$ are coupled in the scheme \eqref{eq:ade_st5} to compute the next step $\discreteFunc{\phi}{i}{j}{n+1}$ of $\discreteFunc{\phi}{i}{j}{n}$.
We substitute the formula \eqref{eq:ade_st4} into \eqref{eq:ade_st3} and obtain the following scheme for $\discreteFunc{\phi}{i}{j}{n+1}$:
\begin{equation*}
    \discreteFunc{\phi}{i}{j}{n+1} = \sum_{\substack{i-2\leq k\leq i+2\\ j-2\leq \ell\leq j+2}}\discreteFunc{C}{k}{\ell}{n}\discreteFunc{\phi}{k}{\ell}{n},
\end{equation*}
where
\begin{equation*}
    \discreteFunc{C}{i}{j}{n} = 1 - \frac{\Delta t}{\Delta x} - \frac{\Delta t}{\Delta y} - \frac{1+\beta_y}{2|B|}\left(\frac{\Delta t}{\Delta x}\right)^2\alpha_s - \frac{1+\beta_x}{2|B|}\left(\frac{\Delta t}{\Delta y}\right)^2\alpha_s,
\end{equation*}
\begin{align}\label{eq:coff}
    \discreteFunc{C}{i+1}{j}{n} = \discreteFunc{C}{i-1}{j}{n} = \frac{\Delta t}{2\Delta x}, &\quad \discreteFunc{C}{i+2}{j}{n} = \discreteFunc{C}{i-2}{j}{n} = \frac{1+\beta_x}{|B|}\left(\frac{\Delta t}{2\Delta x}\right)^2\alpha_s,\nonumber\\
    \discreteFunc{C}{i}{j+1}{n} = \discreteFunc{C}{i}{j-1}{n} = \frac{\Delta t}{2\Delta y}, &\quad \discreteFunc{C}{i}{j+2}{n} = \discreteFunc{C}{i}{j-2}{n} = \frac{1+\beta_y}{|B|}\left(\frac{\Delta t}{2\Delta y}\right)^2\alpha_s,\nonumber\\
    \discreteFunc{C}{i+1}{j+1}{n} = \discreteFunc{C}{i-1}{j-1}{n} = \frac{\beta_c}{2|B|}\frac{\Delta t}{\Delta x}\frac{\Delta t}{\Delta y}\alpha_s, &\quad \discreteFunc{C}{i+1}{j-1}{n} = \discreteFunc{C}{i-1}{j+1}{n} = -\frac{\beta_c}{2|B|}\frac{\Delta t}{\Delta x}\frac{\Delta t}{\Delta y}\alpha_s
\end{align}
and the other coefficients $\discreteFunc{C}{k}{\ell}{n}$ are zero.
\begin{rem}\label{rem:mon}
Due to the cross terms in \eqref{eq:coff}, the scheme is not monotone for any choice of $\alpha_s \neq 0$ for symmetric matrices $\matK$.
However, these terms are of order $O(\alpha_s^2)$ as $\alpha_s\to 0$,
and hence we can expect that the scheme is \textit{approximately monotone for sufficiently small $\alpha_s$}.
On the other hand, the scheme \eqref{eq:ade_st3} suffers from the terms having negative coefficients due to its order $O(\alpha_s)$,
and hence we can say that the scheme \eqref{eq:ade_st5} is refined to be more monotone than the original one.
\end{rem}
\begin{rem}\label{rem:zero}
The assumption that $\discreteFunc{u}{i}{j}{n} = \discreteFunc{v}{i}{j}{n}\equiv 0$ is not really feasible in applications since this implies that
the spatial gradients of $\phi$ equals zero. This assumption is crucial to analyze the DMP of the scheme \eqref{eq:ade_st5} in the next section,
and is not invoked in the numerical experiments (see \Section{sec:num}).
\end{rem}

\section{Discrete maximum principle of FDM}
\subsection{Equivalent condition for DMP}
The three conditions of Barles and Souganidis \cite{BarlesSouganidis1991} are sufficient criteria to guarantee the convergence of a discrete scheme,
and these conditions might not be appropriate to establish the discrete maximum principle (DMP).
In this section, we recall an equivalent condition for DMP from Ciaret \cite{Ciarlet1970}.

Let $\numBulk$ and $\numBdd$ be the numbers of mesh points in the interior and on the boundary, respectively; that is to say 
\begin{equation*}
    \numBulk := N_xN_y - 4(N_x + N_y - 4)\quad\mbox{and}\quad \numBdd := 4(N_x + N_y - 4).
\end{equation*}
We suppose that the first equation of \eqref{eq:ade} is represented as $L\phi = 0$ by using the elliptic differential operator $L$.
Then, we approximate the operator $L$ by the discrete differential operator $L^h$ defined by
\begin{equation}\label{eq:fdm}
    (L^h\phi)_{i} :=  \sum_{j=1}^{\numBulk} a_{ij}\phi_{j} + \sum_{j=1}^{\numBdd}a^\partial_{ij}\phi^\partial_{j}\qquad\mbox{for}\quad 1\leq i\leq N,
\end{equation}
where $\phi_{j}\,(1\leq j\leq \numBulk)$ and $\phi^\partial_{j}\,(1\leq j\leq \numBdd)$ are the values of $\phi$ at the mesh points in the interior and on the boundary of $\Omega$, respectively.
We let $\mat{A} := \left[a_{ij}\right]_{i,j}\in\mathbb{R}^{\numBulk^2}$ and $A^\partial := \left[a^\partial_{ij}\right]_{i,j}\in\mathbb{R}^{\numBulk\times\numBdd}$.
Let us recall the definition of DMP for the discrete differential operator $L^h$ from Ciaret \cite{Ciarlet1970}:

\begin{dfn}[Discrete maximum principle]\label{dfn:dmp}
    The discrete differential operator $L^h$ satisfies DMP provided that 
    if $(L^h\phi)_i\leq 0$ for every $1\leq i\leq N$, then it follows that
    \begin{equation*}
        \max\left\{\phi_j\biggm| 1\leq j\leq N\right\} \leq 0\lor\max\left\{\phi^\partial_j\biggm| 1\leq j\leq N^\partial\right\}.
    \end{equation*}
\end{dfn}

To state the result of Ciaret \cite{Ciarlet1970}, we introduce the following matrices:
\begin{equation*}
    \mat{\widetilde{A}} := \left[
		\begin{array}{cc}
        	A & A^\partial\\
        	O & I^\partial
		\end{array}
	\right]
    \in\mathbb{R}^{(\numBulk + \numBdd)^2}
    \quad\mbox{and}\quad \mat{\widetilde{G}} := \left[
		\begin{array}{cc}
			G & G^\partial\\
			O & I^\partial
		\end{array}
	\right]
    \in\mathbb{R}^{(\numBulk\times\numBdd)^2},
\end{equation*}
where $I^\partial\in\mathbb{R}^{(\numBdd)^2}$ and $O$ respectively denote the identity matrix and the zero matrix;
$G := \mat{A}^{-1}$ and $G^\partial := -\mat{A}^{-1}\mat{A^\partial}$. Then, we easily observe that $\mat{\widetilde{G}} = \mat{\widetilde{A}}^{-1}$.
We are now ready to state an equivalent condition for DMP which was shown by Ciaret \cite{Ciarlet1970}:
\begin{prop}[THEOREM 1 of \cite{Ciarlet1970}]\label{prop:dmp}
    The discrete differential operator $L^h$ satisfies the discrete maximum principle if and only if
    \begin{equation*}
        \mat{\widetilde{G}}\geq\mat{O} \qquad\mbox{and}\qquad \mat{G^\partial}\vec{e}^\partial\leq\vec{e}
    \end{equation*}
    hold, where the second inequality means that all corresponding components of the vectors satisfy the inequality;
    we let $\vec{e}\in\mathbb{R}^{\numBulk}$ and $\vec{e}^{\partial}\in\mathbb{R}^{\numBdd}$ be the vectors whose components are all equal to $1$.
\end{prop}
In what follows, we suppose that meshing of the domain $\Omega$ is minimal; that is to say $N_x = N_y = 3$, then we have $\numBulk = 1$ and $\numBdd = 8$.
For simplicity, we suppose that $h = \Delta x = \Delta y$.

\subsection{Optimal choice of parameters for DMP preservation}
Let us investigate a relationship on the parameters $\alpha_s$, $k_x$, $k_y$, $k_c$, $\Delta t$ and $h$ to satisfy DMP for the hyperbolic-type scheme \eqref{eq:ade_st5}.
First, we observe that the scheme \eqref{eq:ade_st5} is represented as
\begin{equation}\label{eq:ade_st6}
    \discFunc{\phi}{n+1} := \discFunc{\phi}{n} + L^h_0\discFunc{\phi}{n},
\end{equation}
where the discrete differential operator $L^h_0$ is defined by the formula \eqref{eq:fdm} with
$A\in \mathbb{R}^{1\times 1}$ and $A^\partial\in\mathbb{R}^{1\times 8}$ given by
\begin{equation*}
    \mat{A} = \left[-\frac{2\Delta t}{h} - \frac{1+\beta_y}{2|B|}\left(\frac{\Delta t}{h}\right)^2\alpha_s - \frac{1+\beta_x}{2|B|}\left(\frac{\Delta t}{h}\right)^2\alpha_s\right]
\end{equation*}
and
\begin{multline*}
    \mat{A^\partial} = \Bigg[
        \frac{\beta_c}{2|B|}\left(\frac{\Delta t}{h}\right)^2\alpha_s,\ \frac{\Delta t}{2h} + \frac{1+\beta_x}{|B|}\left(\frac{\Delta t}{2h}\right)^2\alpha_s, 
        -\frac{\beta_c}{2|B|}\left(\frac{\Delta t}{h}\right)^2\alpha_s,\ \frac{\Delta t}{2h} + \frac{1+\beta_c}{|B|}\left(\frac{\Delta t}{2h}\right)^2\alpha_s,\\
        \qquad\qquad \frac{\beta_c}{2|B|}\left(\frac{\Delta t}{h}\right)^2\alpha_s,\ \frac{\Delta t}{2h} + \frac{1+\beta_c}{|B|}\left(\frac{\Delta t}{2h}\right)^2\alpha_s,\ 
        -\frac{\beta_c}{2|B|}\left(\frac{\Delta t}{h}\right)^2\alpha_s,\ \frac{\Delta t}{2h} + \frac{1+\beta_y}{|B|}\left(\frac{\Delta t}{2h}\right)^2\alpha_s
    \Bigg].
\end{multline*}
Thus, $\mat{\widetilde{G}}\geq \mat{O}$ if and only if 
\begin{equation}\label{eq:dg}
    \frac{2\Delta t}{h} + \left\{\frac{1+\beta_x}{2|B|}\left(\frac{\Delta t}{h}\right)^2 + \frac{1+\beta_y}{2|B|}\left(\frac{\Delta t}{h}\right)^2\right\}\alpha_s\leq 0.
\end{equation}
Meanwhile, since $\mat{A}\geq \mat{O}$, the condition $\mat{G^\partial}\vec{e}^\partial\leq\vec{e}$ is equivalent to
$-\sum \mat{A^\partial}\leq\mat{A}$, and thus we have
\begin{equation}\label{eq:dg2}
    \frac{\trace{\matK} - 2k_c}{|B|}\alpha_s \leq 0.
\end{equation}
Noting that $\trace{\matK}^2 - (2k_c)^2 = (k_x - k_y)^2 + 4\detK > 0$, the condition \eqref{eq:dg2} is equivalent to saying that
$\alpha_s$ and $|B|$ have the different signs.
In the sequel, we shall prove that there exists $\alpha_s$ that fulfills the conditions \eqref{eq:dg} and \eqref{eq:dg2} under some assumptions on $\Delta t$ and $h$:
\begin{thm}\label{thm:dmp}
    Assume that $\discFunc{\phi}{n+1} = \discFunc{\phi}{n}$ for some $n\in\mathbb{N}$.
    For each $C \in \mathbb{R}$, let $\thresTildeM(C)$ and $\thresTildeP(C)$ be the two roots of the equation $f_C(\alpha) = 0$, where
    \begin{equation}\label{eq:thres2}
        f_C(\alpha) := (4\Delta t^2 + C\trace{\matK}\Delta t)\alpha^2 + (4\trace{\matK}\Delta t + 2C\detK)\alpha + 4\detK.
    \end{equation}
    Let $I_{DMP}(C)\subset\mathbb{R}$ be defined by
    \begin{equation*}
        I_{DMP}(C) := \intvalClosed{\thresTildeM(C)}{\thresM}\cup\intvalClosed{\thresP}{\thresTildeP(C)},
    \end{equation*}
    where $\thresM$ and $\thresP$ have been defined in \eqref{eq:thres}.
    Then, the scheme \eqref{eq:ade_st5} satisfies DMP if and only if $\alpha_s \in I_{DMP}\left(\Delta t/h\right)$.
    Moreover, the set $I_{DMP}\left(\Delta t/h\right)$
    is not empty whenever
    \begin{equation}\label{eq:thres3}
        \Delta t \leq \frac{\trace{\matK}^2 - \trace{\matK}\sqrt{\trace{\matK}^2 - 4\detK}}{2}.
    \end{equation}
\end{thm}
\begin{rem}
    It is easily seen that $I_{DMP}(0) \neq \emptyset$ since $\thresTildeM(0) = \thresM$ and $\thresTildeP(0) = \thresP$, and thus $I_{DMP}(0) = \{\thresP, \thresM\}$.
    Here, we need to show that $I_{DMP}(C)$ is increasing if $C$ is close to $0$.
\end{rem}
\begin{proof}[\textbf{Proof of \Theorem{thm:dmp}}]
    For short notation, we hereafter write $C := \Delta t/h$. First, we can exclude the case when $\alpha_s > 0$ and $|B| < 0$ since the latter condition is equivalent to $\alpha_s \in \intval{\thresM}{\thresP}\subset\intval{-\infty}{0}$.
    Thus, we only address the case when $\alpha_s < 0$ and $|B| > 0$.
    We easily observe that the equation \eqref{eq:thres2} has two roots $\thresTildeM(C)$ and $\thresTildeP(C)$
    for sufficiently small $C$ represented as
    \begin{equation*}
        \thresTildePM(C) := \frac{-(2\trace{\matK}\Delta t + C\detK)\pm \sqrt{D(C)}}{4\Delta t^2 + C\trace{\matK}\Delta t},
    \end{equation*}
    where $D(C)$ denotes the determinant of the quadratic equation \eqref{eq:thres2}:
    \begin{equation*}
        D(C) := (2\trace{\matK}\Delta t + C\detK)^2 - 4\detK(4\Delta t^2 + C\trace{\matK}\Delta t).
    \end{equation*}
    Indeed, we compute
    \begin{multline*}
        D(0) = (2\trace{\matK}\Delta t)^2 - 4\detK\cdot 4\Delta t^2 = 4\trace{\matK}^2\Delta t^2 - 16\detK\Delta t^2\\
        = 4\Delta t^2(\trace{\matK}^2 - 4\detK) = 4\Delta t^2\{(k_x - k_y)^2 + 4k_c^2\} > 0.
    \end{multline*}
    Since the signs of $\alpha_s$ and $|B|$ are different each other, the condition \eqref{eq:dg} is equivalent to saying that
    $\alpha := \alpha_s$ is a subsolution of the equation \eqref{eq:thres2}, say $\thresTildeM(C)\leq \alpha_s \leq \thresTildeP(C)$.
    Meanwhile, the condition that $|B| \geq 0$ is equivalent to either $\alpha_s \leq \thresM$ or $\thresP \leq \alpha_s$.
    Hence, the former assertion has been confirmed.

    To ensure that such $\alpha_s$ exists, we need to check that $I_{DMP}(C)\neq \emptyset$ for sufficiently small $C$.
    For this purpose, we first note that $\thresTildePM(0) = \thresPM$ and show that $\frac{d}{dC}{\thresTildePM}(0) \gtrless 0$.
    Since $\thresTildeM(C)$ and $\thresTildeP(C)$ are smooth, this implies that $I_{DMP}(C)$ is increasing with respect to $C$
    in a neighborhood of $C = 0$, and thus $I_{DMP}(C)\neq \emptyset$ for sufficiently small $C$.

    A direct calculation shows that
    \begin{multline*}
        \frac{d}{dC}\thresTildePM(C) = \frac{\left(-\detK\pm \frac{D'(C)}{2\sqrt{D(C)}}\right)(4\Delta t^2 + C\trace{\matK}\Delta t) + \left\{(2\trace{\matK}\Delta t + C\detK)\mp\sqrt{D(C)}\right\}\trace{\matK}\Delta t}{(4\Delta t^2 + C\trace{\matK}\Delta t)^2}\\
        =\frac{(-2\sqrt{D(C)}\Delta t\pm D'(C))(4\Delta t + C\trace{\matK}) + (2\sqrt{D(C)}(2\trace{\matK}\Delta t + C\detK)\mp 2D(C))\trace{\matK}}{(4\Delta t + C\trace{\matK})^2\Delta t\cdot 2\sqrt{D(C)}}
    \end{multline*}
    Let us analyze the sign of the numerator $\mathcal{N}^\pm(0)$ of $\frac{d}{dC}\thresTildePM(0)$.
    Using the notation that $E(C) := \sqrt{D(C)}$ for $C > 0$, the numerator is represented as
    \begin{equation}\label{eq:num}
        \mathcal{N}^\pm(0) = \mp 2\trace{\matK}E(0)\left\{E(0) + \frac{-2\Delta t(2\Delta t - \trace{\matK}^2)}{\mp\trace{\matK}}\right\}.
    \end{equation}
    Here, we have used that
    \begin{equation*}
        D'(C) = 2(2\trace{\matK}\Delta t + C\detK)\cdot \detK - 4\detK\trace{\matK}\Delta t,
    \end{equation*}
    and hence $D'(0) = 0$.
    Let $E^\pm(0)$ be the two roots of the equation $\mathcal{N}^\pm(0) = 0$. Then, we have
    \begin{equation*}
        E^+(0) = \frac{2\Delta t(\trace{\matK}^2 - 2\Delta t)}{\trace{\matK}}.
    \end{equation*}
    We observe that for $\Delta t < \trace{\matK}^2/2$, the inequality $E(0)\leq E^+(0)$
    is equivalent to saying that
    \begin{equation*}
        \Delta t^2 - \trace{\matK}^2\Delta t + \trace{\matK}^2\detK \geq 0.
    \end{equation*}
    It is easily seen that the above inequality is satisfied if either $\Delta t \leq \Delta t^-$ or $\Delta t^+\leq \Delta t$, where
    \begin{equation*}
        \Delta t^- := \frac{\trace{\matK}^2 - \trace{\matK}\sqrt{\trace{\matK}^2 - 4\detK}}{2},\quad\mbox{and}\quad \Delta t^+ := \frac{\trace{\matK}^2 + \trace{\matK}\sqrt{\trace{\matK}^2 - 4\detK}}{2}.
    \end{equation*}
    Since we have supposed that $\Delta t < \trace{\matK}^2/2$, if $\Delta t\leq \Delta t^-$, then $E(0)\leq E^+(0)$.
    For such $\Delta t$, we see that $\mathcal{N}^+(0) > 0$, and hence it follows that $\frac{d}{dC}\thresTildeP(0) > 0$.
    Meanwhile, if $\Delta t\leq \Delta t^-\leq \trace{\matK}^2/2$, the inequality $\mathcal{N}^-(0) < 0$ is straightforward, and we deduce that $\frac{d}{dC}\thresTildeM(0) < 0$.
    This completes the proof.
\end{proof}
\begin{rem}
Replacing $\alpha_s$ with $-\alpha_s$ in the above arguments, we may assume that $\alpha_s > 0$ without loss of generality.
Hence, we hereafter suppose that $\alpha_s > 0$. In this case, $\thresTildeP$ and $\thresTildeM$ should exchange their roles,
although they will be written in the same manner for the reader's convenience.
\end{rem}
\begin{exm}\label{exm:thres1}
    Let us calculate a threshold value of $\alpha_s$ to ensure DMP for the scheme \eqref{eq:ade_st5} (see again \Eqref{eq:mobility3}).
    Here, we revisit the example in \Section{sec:basic} (see \Eqref{eq:aniso2}) when the anisotropic tensor is given by
    \begin{equation*}
        \matK = \left[
			\begin{array}{cc}
				k_x & k_c\\
				k_c & k_y
			\end{array}
		\right]
         =
        \Theta^{-1} \left[
			\begin{array}{cc}
            	1 & 0\\
            	0 & \left(\frac{k_{||}}{k_\perp}\right)^{-1}
			\end{array}
		\right]
        \Theta\qquad\mbox{with}\quad \Theta := \left[
			\begin{array}{cc}
            	\cos{\theta} & -\sin{\theta}\\
            	\sin{\theta} & \cos{\theta}
			\end{array}
		\right].
    \end{equation*}
    Namely, we shall treat the case when the magnetic field lines lie horizontally in the domain $\Omega$ at an angle $\theta$.
    In particular, we set $\theta := \pi / 4$ and $k_{||}/k_\perp := 10^4$. Then, we see that $k_x=k_y=0.50005$ and $k_c = 0.49995$.
    We stress that the time step $\Delta t$ is set according to the condition \eqref{eq:thres3}. Then, we have $\Delta t\approx 0.0001$.
    According to the discussion in \Section{sec:mon}, the smallness of $|\alpha_s|$ is important to ensure the convergence of the scheme \eqref{eq:ade_st5}.
    Thus, we only consider $\thresP$ and $\thresTildeP(C)$ since these are much closer to zero than $\thresM$ and $\thresTildeM(C)$, respectively.

    Fixing $\Delta t := 0.0001$, we show optimal parameters $\thresTildeP(C)$ for several settings of $C$ and $h$ in \Table{tab:thres}.
	Here, $h$ is determined according to the CFL condition, say $C = \Delta t/ h$ for each $C > 0$.

    \begin{table}[H]
        \centering
        \begin{tabularx}{0.6\textwidth}{|X|X|X|} 
            \hline
            \textbf{C} & \textbf{h} & \textbf{$\thresTildeP(C)$} \\
            \hline\hline
            0.5      & 0.0002      & 0.8768     \\
            0.25      & 0.0004      & 0.9376     \\
            0.1      & 0.001      & 0.9749     \\
            0.05      & 0.002      & 0.9874     \\
            0.025      & 0.004     & 0.9936     \\
            0.01      & 0.01      & 0.9974     \\
            \hline
        \end{tabularx}
        \caption{Optimal values of $\thresTildeP(C)$ for DMP}
        \label{tab:thres}
    \end{table}
\end{exm}
\begin{rem}
A direct calculation shows that $\thresTildeP(C)\to 0$ (resp. $\thresTildeP(C) \to 1$) as $C\to \infty$ (resp. $C\to 0$).
\end{rem}

	\section{Numerical tests}\label{sec:num}
	The analysis in the previous section proves that the scheme parameter $\alpha_{\rm s}$ arising in the hyperbolic system approach is related to the DMP preservation.
	It is also implied that, with an appropriate choice of $\alpha_{\rm s}$, one can obtain DMP-preserving numerical solutions even without a nonlinear discretization treatment, such as a priori or a posteriori limiters.
	It is noted that the DMP analysis in the previous section is based on \Eqref{eq:ade_st5}, where the contributions from $u^{n}$ and $v^{n}$ (i.e., information from the previous-step subsequent variables') are neglected.
	Accordingly, the value of $\alpha_{\rm s}$ that preserves the DMP in practice is expected to be close to the analysis, though there be a discrepancy between the actual numerical simulation and the mathematical analysis.
	Therefore, we verify that the optimal $\alpha_{\rm s}$ exists near the analytical threshold values such as those in Table \ref{tab:thres} through several numerical test cases. 
	
\begin{figure}[t]
		\begin{center}
			\includegraphics[width=160mm]{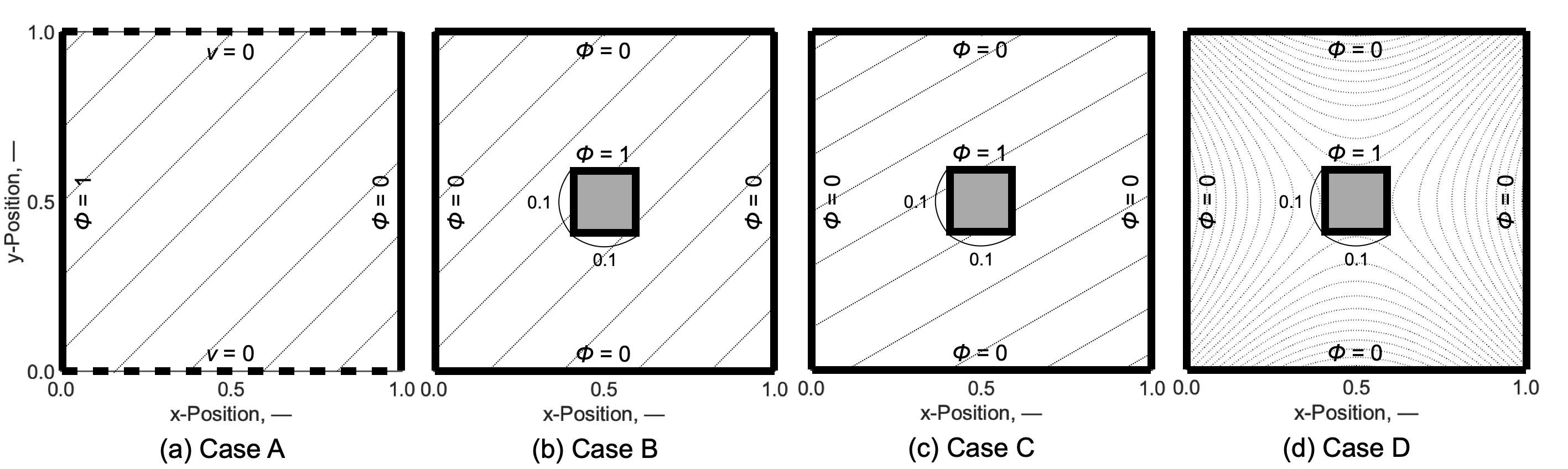}
		\end{center}
		\caption{Two-dimensional test case conditions. Case A and B assume $\pi/4$ anisotropy (magnetic field lines), Case C assumes $\pi/6$, and Case D employs a cusped magnetic field geometry.}
		\label{fig:condition}
	\end{figure}

	The calculation conditions for the two-dimensional test problems are shown in \Figure{fig:condition}.
	Case A corresponds to the DC discharge plasma devices with magnetic confinement, and similar setups appear in Refs. \cite{kawashima2015,mikellides2011}.
	A potential difference is imposed between the left and right boundaries.
	The top and bottom boundaries are impermeable, which is implemented by the Dirichlet condition $v=0$ in the hyperbolic system, since $v$ corresponds to the y-direction flux.
	Case B is also a benchmark assessing DMP preservation of schemes, and similar problems are used in Refs. \cite{kuzmin2009,liska2008}.
	A potential difference is applied between the central square region and the outer boundaries.
	Case C uses the same boundary conditions as Case B, but inclines the magnetic field lines by $\pi/6$.
	Case D assumes a cusped magnetic field geometry, with magnetic vector potential
	\begin{equation}
		A_z = \frac{k_{||}}{k_\perp}\left(\left(x-0.5\right)^2-\left(y-0.5\right)^2\right).
	\end{equation}
	The magnetic field is given by $\vec{B} = \nabla \times A_z$.
	This configuration also produces a steep distribution of $\phi$, making DMP preservation challenging.
	The strength of the anisotropy, $k_{||}/k_{\perp}$, is set to 10$^4$ in all cases.
	For each test, a numerical solution is verified to be DMP preserving if $0\leq\phi\leq 1$ everywhere in the simulation domain.

	All test cases are computed using the preconditioned hyperbolic finite-difference scheme.
	The grid sizing is $h=0.01$ and pseudo-time step is $\Delta t = 10^{-4}$.
	According to the analysis result summarized in Table 1, the optimal $\alpha_{\rm s}$ is 0.9974.
	This result suggests that DMP-preserving solutions should be obtained when $\alpha_{\rm s}$ is larger than 0.9974.
	Thus, we vary $\alpha_{\rm s}$ from 0.1 to 4.0, and calculated the test cases.
	\begin{figure}[t]
		\begin{center}
			\includegraphics[width=125mm]{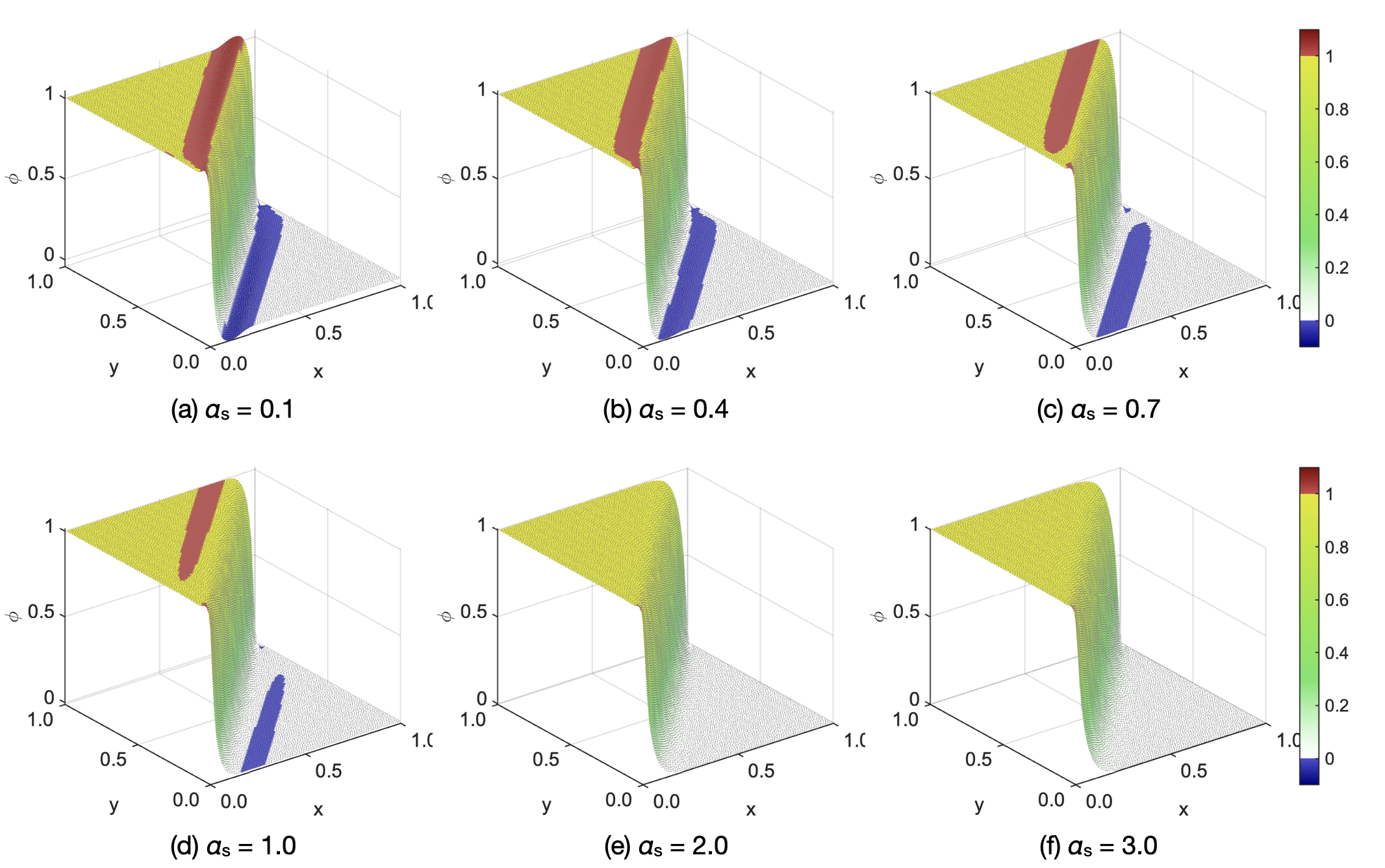}
		\end{center}
		\caption{Calculation results of anisotropic diffusion test case A computed with the hyperbolic system approach with varied $\alpha_{\rm s}$.}
		\label{fig:phi3d_a}
		\begin{center}
			\includegraphics[width=110mm]{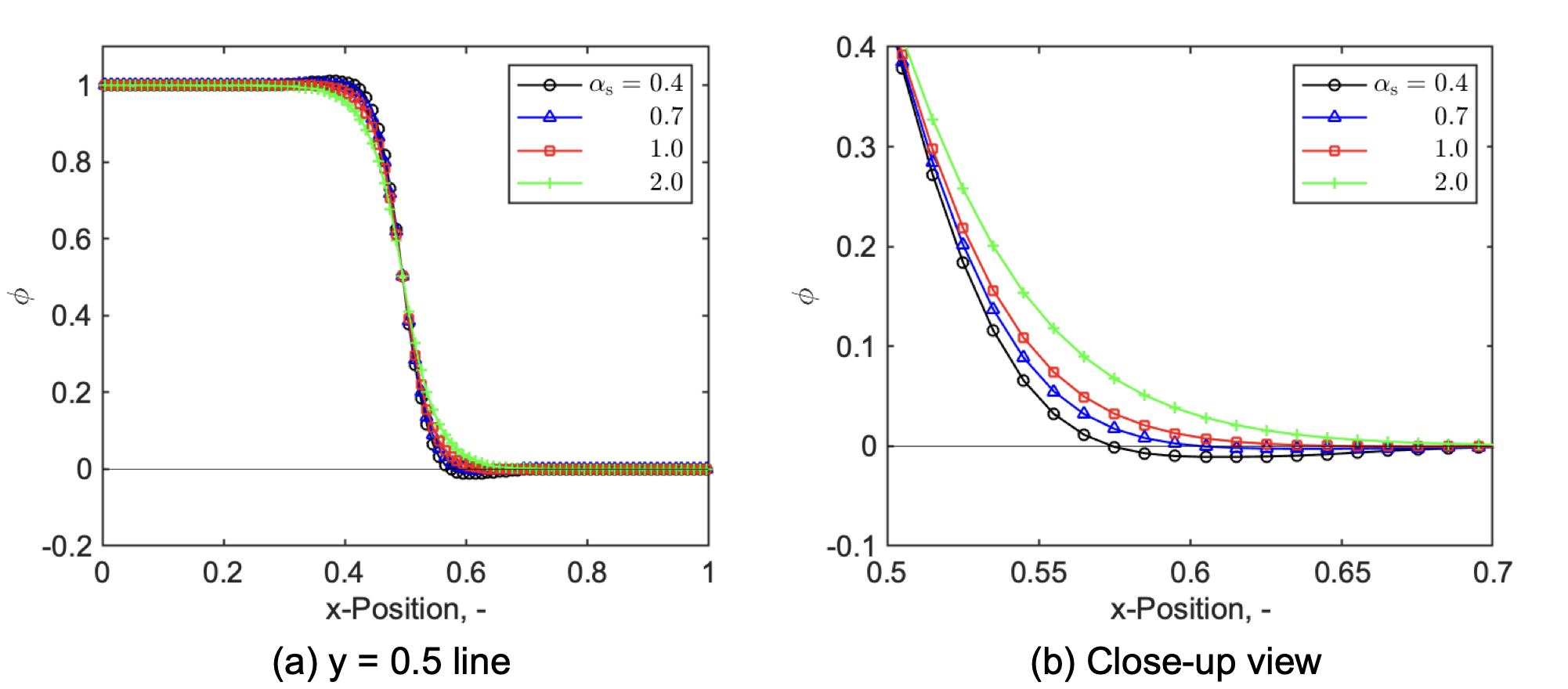}
		\end{center}
		\caption{$\phi$ distribution along $y=0.5$ for Case A in \Figure{fig:condition}, with varying $\alpha_{\rm s}$.}
		\label{fig:phi1d_a}
	\end{figure}

	\Figure{fig:phi3d_a} shows the Case A results as $\alpha_{\rm s}$ is varied.
	Red and blue regions represent overshoot ($\phi > 1$) and undershooting ($\phi < 0$) areas, respectively.
	For smaller $\alpha_{\rm s}$, the results exhibit relatively sharp potential drops accompanied by overshoot and undershoot, indicating DMP violation.
	Even at $\alpha_{\rm s} = 1$, regions where $\phi > 1$ and $\phi < 0$ remain.
	As $\alpha_{\rm s}$ increases further, DMP-preserving results are achieved as shown in \Figure{fig:phi3d_a}(e) and (f).

	The one-dimensional profiles of $\phi$ along $y = 0.5$ line are shown in \Figure{fig:phi1d_a},
	with a close-up view in \Figure{fig:phi1d_a}(b), highlighting the extent of undershoot ($\phi < 0$).
	Undershoots are clearly mitigated as $\alpha_{\rm s}$ increases.
	The value $\alpha_{\rm s}=1.0$ appears close to the threshold between DMP preservation and violation.
	The analysis result in \Section{sec:mon} suggests that the optimal $\alpha_{\rm s}$ is 0.9974.
	A slight discrepancy in $\alpha_{\rm s}$ that preserves DMP is thus observed between the analytical prediction and the numerical results.
	
	\begin{figure}[t]
		\begin{center}
			\includegraphics[width=125mm]{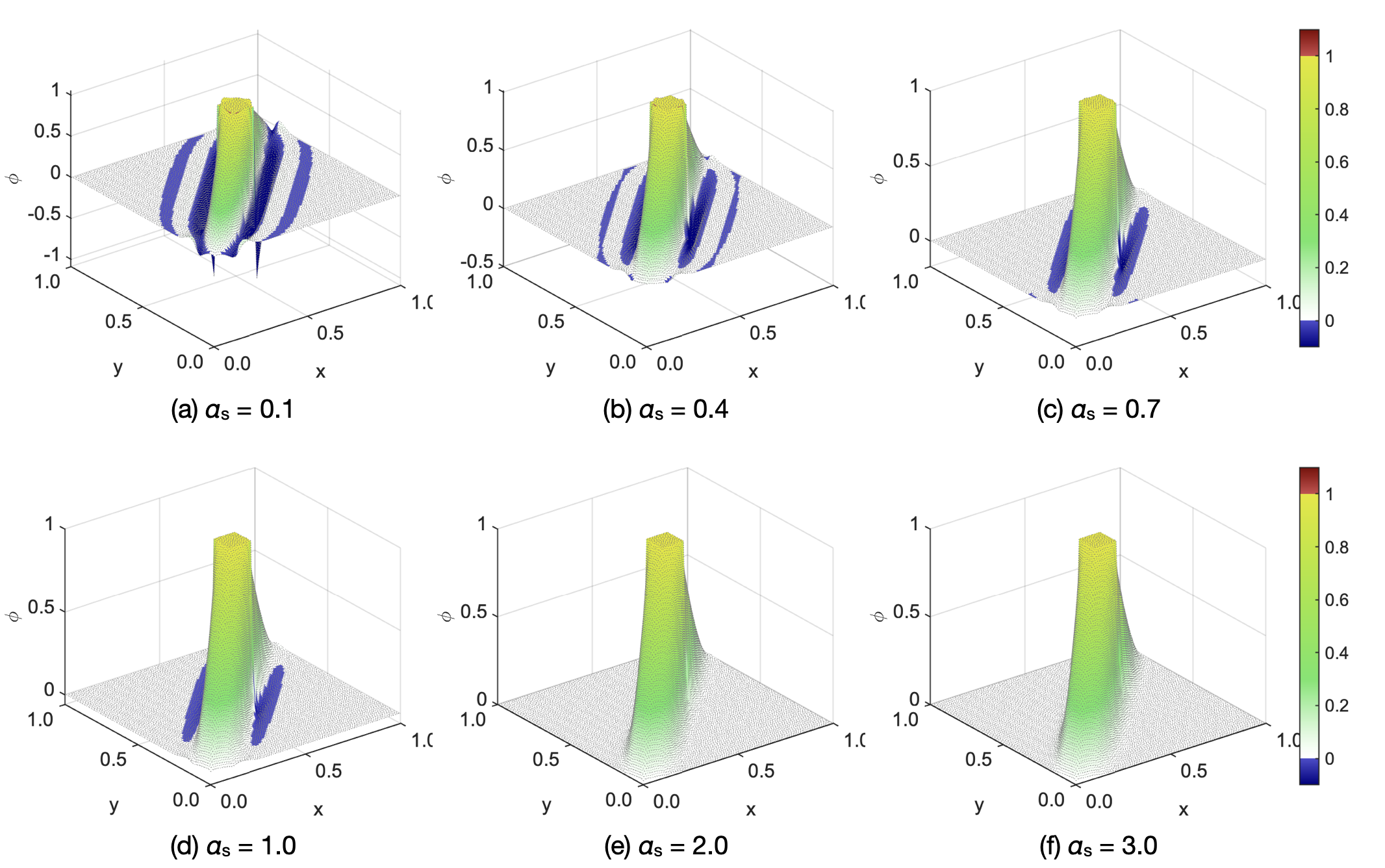}
		\end{center}
		\caption{Calculation results of anisotropic diffusion test case B solved by the hyperbolic system approach with varied $\alpha_{\rm s}$.}
		\label{fig:phi3d_b}
		\begin{center}
			\includegraphics[width=110mm]{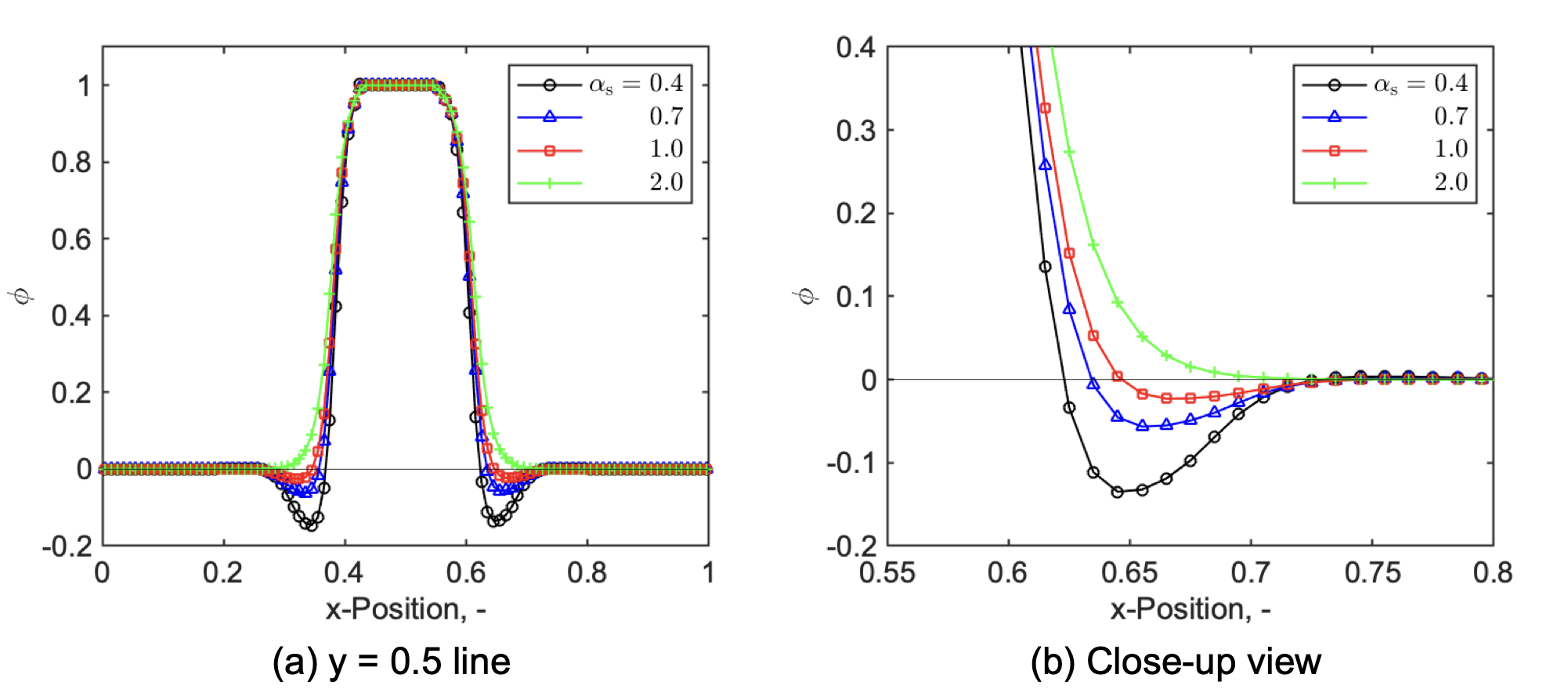}
		\end{center}
		\caption{$\phi$ distribution along $y=0.5$ for Case B in \Figure{fig:condition}, with varying $\alpha_{\rm s}$.}
		\label{fig:phi1d_b}
	\end{figure}

	The calculation results for Case B are presented in \Figure{fig:phi3d_b}.
	In this case, the blue region ($\phi<0$) serves as an indicator of DMP violation.
	As in Case A, negative $\phi$ occurs when $\alpha_{\rm s}$ is 1.0 or smaller, indicating a failure in DMP, whereas for $\alpha_{\rm s}$ exceeding 2.0, the DMP is satisfied.
	The one-dimensional profiles along $y=0.5$ line are shown in \Figure{fig:phi1d_b} (and close-up view in panel (b)).
	The undershoot region where $\phi<0$ shrinks as $\alpha_{\rm s}$ increases, and a DMP-preserving solution is achieved when $\alpha_{\rm s}$ reaches 2.0.

\begin{figure}[t]
		\begin{center}
			\includegraphics[width=125mm]{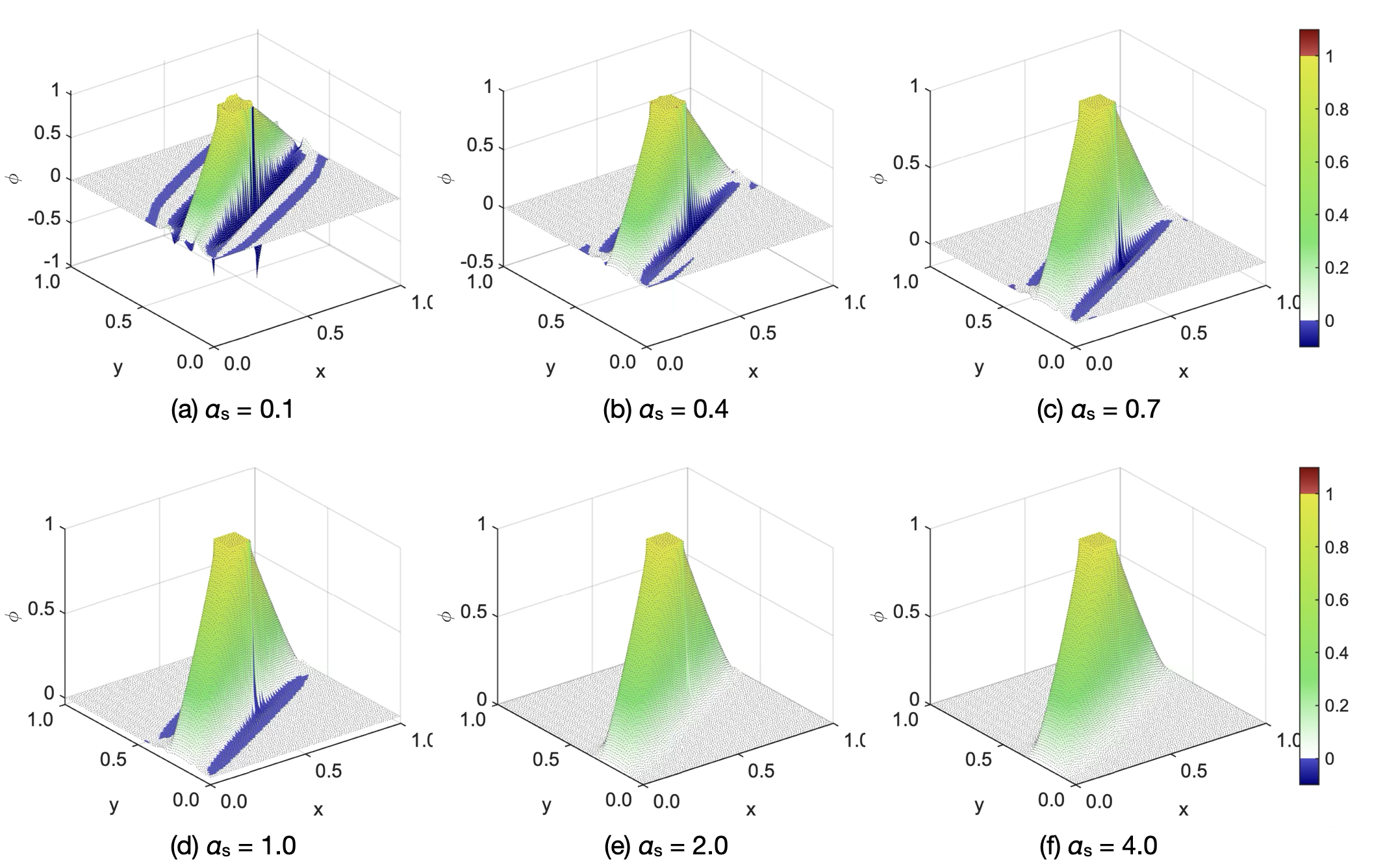}
		\end{center}
		\caption{Calculation results of anisotropic diffusion test case C solved by the hyperbolic system approach with varied $\alpha_{\rm s}$.}
		\label{fig:phi3d_c}
	\end{figure}

\begin{figure}[t]
		\begin{center}
			\includegraphics[width=125mm]{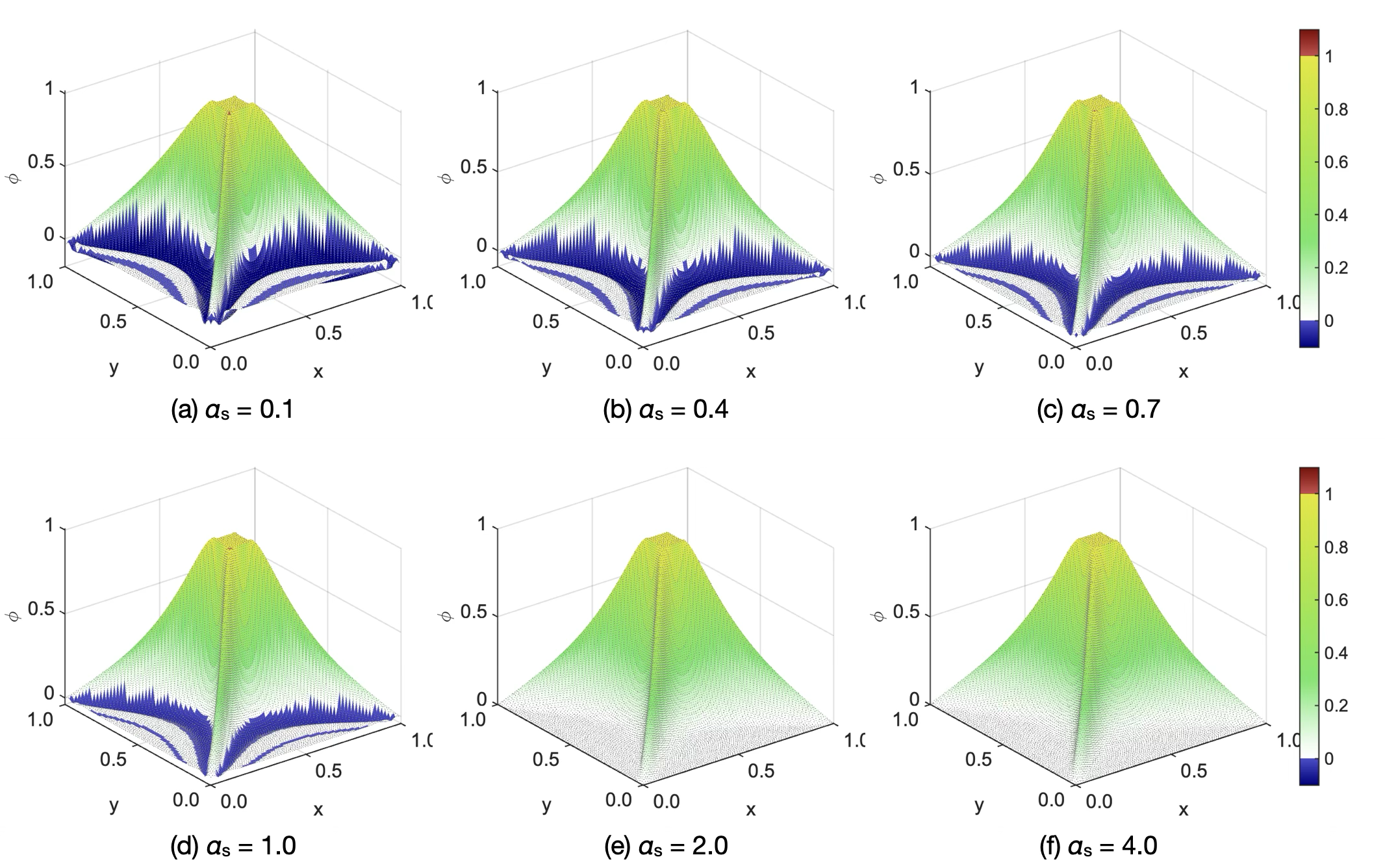}
		\end{center}
		\caption{Calculation results of anisotropic diffusion test case D solved by the hyperbolic system approach with varied $\alpha_{\rm s}$.}
		\label{fig:phi3d_d}
	\end{figure}

	The results of Case C and D are shown in \Figure{fig:phi3d_c} and \Figure{fig:phi3d_d}, respectively.
	In Case C, spurious oscillations are observed when $\alpha_{\rm s}<1$.
	The pattern of undershoot resembles those obtained by elliptic-equation-based schemes without nonlinear limiting \cite{liska2008}. 
	In these cases, the threshold value of $\alpha_{\rm s}$ that preserves the DMP is around 2.0, consistent with Cases A and B.
	This demonstrates the reproducibility of the $\alpha_{\rm s}$ threshold that preserves DMP across different calculation conditions.

\begin{figure}[t]
		\begin{center}
			\includegraphics[width=135mm]{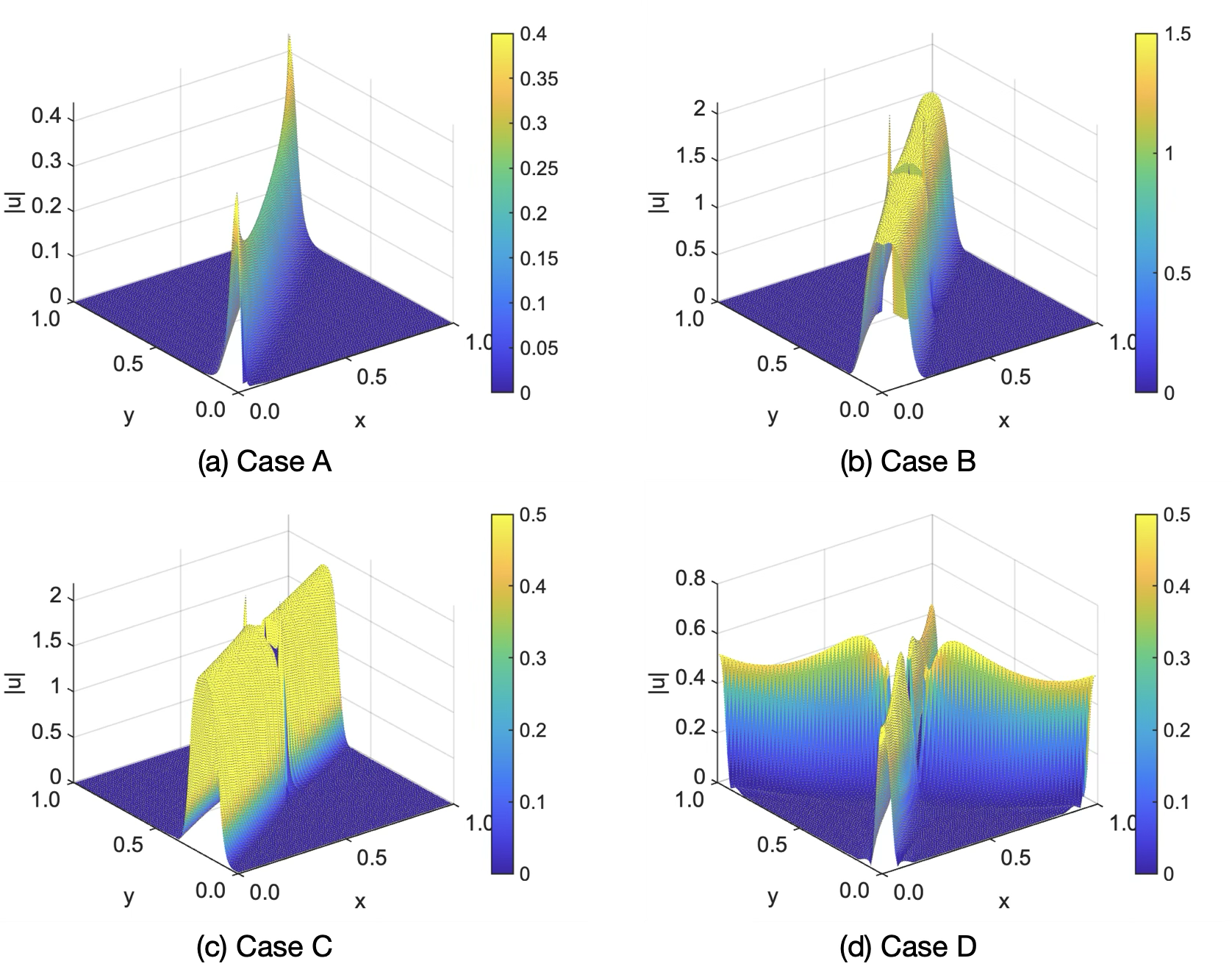}
		\end{center}
		\caption{Distribution of the absolute value of the flow velocity $|u|$ in the test cases A to D obtained with the hyperbolic scheme at $\alpha_{\rm s}=1.0$.}
		\label{fig:uabs3d}
	\end{figure}

	Across these test cases, DMP-preserving solutions are obtained when $\alpha_{\rm s}$ is larger than 2.0.
	This empirical threshold exceeds the analytical value of 0.9974 obtained in \Section{sec:mon}, because the analysis assumes that the subsequent variables of flow velocities $u^n$ and $v^n$ vanish.
	\Figure{fig:uabs3d} shows the absolute velocity $|u|\equiv \sqrt{u^2 + v^2}$ for the four test cases with $\alpha_{\rm s} = 1.0$.
	Since $u$ and $v$ are proportional to the flux (and hence to gradients of $\phi$), they exhibit steep variations where $\phi$ has sharp gradients.
	These steep variations in the flow velocity affects the criterion of DMP, and the larger $\alpha_{\rm s}$ is required.
	Nevertheless, the analysis in \Section{sec:mon} remains useful to determine the range of $\alpha_{\rm s}$ that preserves the DMP.
	The order of magnitude of the optimal $\alpha_{\rm s}$ from analysis is consistent with the numerically obtained threshold.

	In summary, the preconditioned hyperbolic system approach can satisfy the DMP provided that $\alpha_{\rm s}$ is chosen sufficiently large.
	The resulting hyperbolic finite difference scheme remains linear; no a priori flux limiter, a posteriori flux correction, or adaptive stencil selection is required.
	Thus, the hyperbolic system approach offers a simple method to construct DMP-preserving schemes for anisotropic diffusion.
	In this study, $\alpha_{\rm s}$ has been assessed through dimensionless analysis and test problems. 
	For realistic problems with physical units such as the magnetized plasma flows, it is recommended that one applies nondimensionalization to the governing equations, and choose $\alpha_{\rm s}$ near 1.0.

\section{Concluding remarks}
	Anisotropic diffusion equations, which commonly arise in magnetized plasma flows, pose numerical challenges when solved using standard central difference schemes, as these methods often violate the discrete maximum principle (DMP) due to the presence of cross-diffusion terms. 
	Such violation can produce spurious oscillations and nonphysical behavior in numerical solutions, which is particularly important in plasma simulations where accurate potential distributions are essential for predicting charged particle motion. 
	The hyperbolic system approach, which converts the second-order anisotropic diffusion equation into a system of first-order partial differential equations, offers advantages in extensibility and simplicity. 
	However, the properties of this approach with respect to DMP for anisotropic diffusion has not been understood.
	This study aimed to assess the DMP property of the hyperbolic system approach and establish criteria for ensuring monotonicity in anisotropic diffusion problems through mathematical analysis and numerical verification.

	The key contribution of this study is the theoretical framework for determining optimal free parameter so as to preserve the DMP. 
	Through a mathematical analysis, we derived conditions on the artificial parameter $\alpha_{\rm s}$ that appears in the hyperbolic system construction. 
	The analysis revealed that there exists a specific range of $\alpha_{\rm s}$ values that guarantees DMP satisfaction, with the optimal value depending on the diffusion coefficient matrix, time step, and grid spacing. 
	Because the analysis is based on a simplified scheme for the hyperbolic system, numerical tests are performed to verify that the optimal $\alpha_{\rm s}$ exists around the analytically obtained value.
	Numerical experiments across four different test cases confirmed that DMP preservation is achieved when $\alpha_s$ exceeds a threshold value (approximately 2.0 in the present tests), eliminating spurious oscillations and ensuring physically meaningful solutions.

	The practical significance of this work lies in its simplicity and applicability. 
	Unlike other DMP-preserving methods that require nonlinear treatments such as flux limiters or adaptive stencil selection, the proposed hyperbolic formulation with preconditioning and the recommended parameter selection remains a linear scheme. 
	This method retains computational efficiency while providing robust computations for anisotropic diffusion problems. 
	The approach is particularly valuable for plasma simulations where magnetic confinement induces strong anisotropies, and where accurate cross-field transport analysis is critical for predicting device performance. 

\section*{Acknowledgments}
	The authors are grateful to Prof. Daisuke Tagami (Kyushu University) and Dr. Masaru Miyashita for fruitful discussions.
	This work was supported by JSPS KAKENHI Grant Number JP24K01084, and Grant for Short-term Joint Research, Institute of Mathematics for Industry, Kyushu University (No. 2023a022 and 2024a040).
	
\section*{CRediT authorship contribution statement}
	{\bf Tokuhiro Eto:} Writing -- review \& editing, Writing -- original draft, Methodology, Formal analysis; {\bf Rei Kawashima:} Writing -- review \& editing, Writing -- original draft, Methodology, Visualization, Software, Funding acquisition.
	
\section*{Data availability}
	Data will be made available on reasonable request.

\section*{Declaration of competing interest}
	The authors declare that they have no known competing financial interests or personal relationships that could have appeared to influence the work reported in this paper.

\renewcommand{\refname}{Reference}
\bibliographystyle{elsarticle-num}
\bibliography{CMA_Kawashima}

\end{document}